\documentclass[]{article}
\usepackage{graphicx}
\usepackage{graphics}
\usepackage{amssymb,amsmath,amsthm,epsfig}
\newtheorem{theorem}{Theorem}
\newtheorem{lemma}{Lemma}
\newtheorem{corollary}{Corollary}
\newtheorem{definition}{Definition}
\newtheorem{proposition}{Proposition}

\newtheorem{remark}{Remark}
\newcommand{\R}{\mathbb R}

\newcounter{rea}
    \newcommand{\K}{\mathbf K}
\newcommand{\E}{\mathbf E}
\newcommand{\epsn}{\varepsilon_n}

\newcounter{rek}
\begin{document}


%
\begin{center}
{\large {\bf Spectral decay of the  sinc kernel operator and approximations by
Prolate Spheroidal Wave Functions.}}\\
\vskip 1cm Aline Bonami$^a$ and Abderrazek Karoui$^b$ {\footnote{Corresponding author,\\
This work was supported in part by the  ANR grant "AHPI" ANR-07-
BLAN-0247-01, the French-Tunisian  CMCU 10G 1503 project and the
DGRST  research grant 05UR 15-02.}}
\end{center}
\vskip 0.5cm {\small
$^a$ F\'ed\'eration Denis-Poisson, MAPMO-UMR 6628, Department of Mathematics, University of Orl\'eans, 45067 Orl\'eans cedex 2, France.\\
\noindent $^b$ University of Carthage,
Department of Mathematics, Faculty of Sciences of Bizerte, Tunisia.}\\
Emails: aline.bonami@univ-orleans.fr ( A. Bonami), abderrazek.karoui@fsb.rnu.tn (A. Karoui)\\

\noindent{\bf Abstract}--- For fixed $c,$ the  Prolate Spheroidal Wave
Functions (PSWFs) $\psi_{n, c}$ form a basis with remarkable properties
for the space of band-limited functions with bandwidth $c$. They have
been largely studied and used after the seminal work of D. Slepian, H. Landau and H. Pollack.
Recently, they have been used for the approximation of functions
in the  Sobolev space $H^s([-1,1])$. In view of this, we give new estimates on the decay rate
 of eigenvalues of the Sinc kernel integral operators.
This is one of the main issues of this work. A second one  is the choice of the parameter $c$
when approximating a function in  $H^s([-1,1])$ by its truncated PSWFs series expansion.
Such functions may be seen as
the restriction to $[-1,1]$ of almost time-limited  and
band-limited functions, for which PSWFs expansions are still well
adapted.
Finally, we provide the reader with some numerical examples  that
illustrate the different results of this work.\\


\section{Introduction}

For a given value $c>0$, called the bandwidth, PSWFs $(\psi_{n,c}(\cdot))_{n\geq 0}$ constitute an orthonormal basis of $L^2([-1, +1]),$ an orthogonal system of $L^2({\bf R})$ and an
orthogonal  basis of  the Paley-Wiener space $B_c,$ given by $B_c=\left\{ f\in L^2({\bf
R}),\,\, \mbox{Support\ \ } \widehat f\subset [-c,c]\right\}.$
Here, $\widehat f$ denotes the Fourier transform of $f$. They are eigenfunctions of the
  compact integral operators
 $\mathcal F_c$ and $\mathcal Q_c= \frac c{2\pi}\mathcal F_c^*\mathcal F_c$, defined  on $L^2([-1,1])$ by
\begin{equation}\label{eq1.1}
 \mathcal F_c(f)(x)= \int_{-1}^1
e^{i\, c\, x\, y} f(y)\, dy,\quad \mathcal Q_c(f)(x)=\int_{-1}^1\frac{\sin c(x-y)}{\pi (x-y)}\, f(y)\,
dy. \end{equation}
Since   the
operator $\mathcal F_c$ commutes with the Sturm-Liouville operator $\mathcal L_c$,
\begin{equation}\label{eq1.0}
\mathcal L_c(\psi)=-\frac{d}{d\, x}\left[(1-x^2)\frac{d\psi}{d\,x}\right]+c^2 x^2\psi,
\end{equation}
 PSWFs $(\psi_{n,c}(\cdot))_{n\geq 0}$ are also eigenfunctions of $\mathcal L_c$. They are ordered in such a way that the corresponding eigenvalues of $\mathcal L_c$, called $\chi_n(c)$, are strictly increasing. Functions $\psi_{n,c}$ are restrictions to the interval $[-1, +1]$ of real analytic functions on the whole real line and eigenvalues $\chi_n(c)$ are values of  $\lambda$ such that the equation $\mathcal L_c \psi=\lambda \psi$ has a bounded solution.

PSWFs have been introduced  by D. Slepian, H. Landau  and  H. Pollak
\cite{LW, Slepian1, Slepian3, Slepian2} in relation with signal processing. 
 For a detailed review on the  properties, the numerical computations, asymptotic results and
first applications of the PSWFs, the reader is refereed to the recent books on the subject, \cite{Logan}, \cite{Osipov3}.

By Plancherel identity, PSWFs are normalized so that 
\begin{equation}\label{eqq1.4}
\int_{-1}^1 |\psi_{n,c}(x)|^2\, dx = 1,\quad \int_{\mathbb R}
|\psi_{n,c}(x)|^2\, dx =\frac{1}{\lambda_n(c)},\quad n\geq 0.
\end{equation}
Here, $(\lambda_n(c))_n$ is the infinite sequence of  the
eigenvalues of $\mathcal Q_c,$ also arranged in the decreasing order $1>
\lambda_0(c)> \lambda_1(c)>\cdots>\lambda_n(c)>\cdots.$ We call $\mu_n(c)$ the eigenvalues of $\mathcal F_c$. They are given by
$$\mu_{n}(c)=i^n\sqrt{\frac {2\pi}c \lambda_n(c)}.$$
Also, we will adopt the following sign normalization of the PSWFs, given by 
\begin{equation}\label{eeqq1.4}
\psi_{n,c}(0) > 0\mbox{\ \  for even\ \ } n,\quad \, \psi'_{n,c}(0) > 0,\mbox{\ \ for odd \ \ }  n.
\end{equation}

\medskip

One of the main issues that we discuss here is the decay rate of the eigenvalues $\lambda_n(c)$. This decay rate plays a crucial role in most of the 
various concrete applications of the PSWFs.   In this direction, one knows their asymptotic behaviour for $c$ fixed, which has been  given in 1964 by Widom, see \cite{Widom}.
\begin{equation}
 \lambda_n(c)\sim  \left ( \frac {e c} {4(n+\frac 12)} \right)^{2n+1}.
 \end{equation}
This gives the exact decay for $n$ large enough, but one would like to have a more precise information in terms of uniformity of this behaviour, both in $n$ and $c$. On the other hand, Landau has considered  the value of the smallest integer $n$ such that $\lambda_n(c)\leq 1/2$ in \cite{Landau}. More precisely, if we note $c^*_n$ the unique value of $c$ such that $\lambda_n(c)=1/2$, then he proves that
 \begin{equation}
 \label{validityC0intro}
 \frac{\pi}{2}(n-1)\leq c^*_n\leq \frac{\pi}{2} (n+1)\quad\quad \lambda_n(c^*_n)=\frac{1}{2}.
 \end{equation}
So, for $c$ fixed, we almost know when $\lambda_n(c)$ passes through the value $1/2$.  Landau and Widom have also described the asymptotic behaviour, when $c$ tends to $\infty$, of the distribution of the eigenvalues $\lambda_n(c)$.

The search for more precise estimates has attracted a considerable interest, both in numerical and theoretical studies.  We try here to give approximate values for $\lambda_n(c)$ for $c\leq c^*_n$, with some uniformity in the quality of approximation. We rely on the exact formula
\begin{equation}
 \label{eqq2lambda_nintro}
 \lambda_n(c)=\frac{1}{2} \exp\left(-2\int_c^{c^*_n} \frac{(\psi_{n,\tau}(1))^2}{\tau} \, d\tau \right).
  \end{equation}
We use our recent work \cite{Bonami-Karoui1, Bonami-Karoui2}  to estimate the value $\psi_{n,\tau}(1)$. In the first paper it is proved that $|\psi_{n,\tau}(1)|\leq 2\chi_n(\tau)^{1/4}$, which is not sufficient to find a sharp estimate for all values $c$.  The approximation given in the second paper leads to a second estimate, valid for $\frac{\pi n}{2}-c$ larger than some multiple of $\ln n$. We finally find an explicit expression $\widetilde { \lambda_n(c)}$, and prove that it is comparable with $ \lambda_n(c)$ up to some power of $n$. This is given by
\begin{equation}\label{decay}
\widetilde { \lambda_n(c)} = \frac{1}{2} \exp\left(-\frac{\pi^2(n+\frac 12)}{2} \int_{\Phi\left(\frac{2c}{\pi (n+\frac 12)}\right)}^1 \frac{1}{t(\mathbf E(t))^2}\, dt \right)
\end{equation}
Here $\mathbf E$ is the elliptic Legendre integral of the second kind.  The function $\Phi$ is the inverse of the function $t\mapsto \frac{t}{\E(t)}$.

When $n$ tends to $\infty$ with $c$ fixed, we recover the asymptotic behavior given by Widom, which is already a good test of validity. Numerical experiments  prove that this approximation is surprisingly accurate.

As a corollary we have  the following, which may be seen as a kind of quantitative Widom's Theorem.
\begin{theorem} \label{th-intro}
Let $m>0$ be a positive real number and let $M>m,$  $\varepsilon>0$ be  given.  Then there exists a constant $A(\varepsilon, m, M)$ such that, for all $m\leq c\leq M\sqrt n$ and all $n$, we have the inequality
\begin{equation}
 \lambda_n(c)\leq A(\varepsilon, m, M) e^{\varepsilon n}\left( \frac {e c} {4(n+\frac 12)} \right)^{2n+1}.
\end{equation}
\end{theorem}
 One can give an explicit constant $A(\varepsilon, m, M)$. When $c$ may take larger values, there is another statement, where the equivalent found by Widom is replaced by $\widetilde { \lambda_n(c)}$.
 \smallskip

Let us mention that another method to approximate the values $\lambda_n(c)$ has been used by Osipov in \cite{Osipov2}. The estimates given in his paper are of different nature and do not propose such a simple formula. In addition, he mainly considers values of $n$ such that $\frac{\pi n}{2}-c$ is smaller than some multiple of $\ln c$. At this moment both works may be seen complementary. But we underline the fact that numerical tests validate the accuracy of the approximant \eqref{decay}   even when $c$ is close to the critical value, while our theoretical approach is not yet sufficient to do it.
\medskip

Our second
 contribution  is related to the quality of approximation in Sobolev spaces  when a function is replaced by the partial sum of its expansion in some PSWF basis.   This question has attracted a growing interest while, at the same time, were built
PSWFs based numerical schemes for solving various problems from
numerical analysis, see \cite{Boyd1, Boyd2, Chen, Wang, Wang2}.
In particular, in \cite{Boyd1}, the author has shown that a PSWF
approximation based method outperforms in terms of spatial
resolution and stability of timestep, the classical approximation
methods based on Legendre or Tchebyshev polynomials. The authors
of \cite{Chen} were among the first to study the quality of
approximation by the PSWFs in the Sobolev space $H^s(I),
s>0, \, I=[-1,1].$ In particular, they have given an estimate of the decay of
the PSWFs expansion coefficients of a function $f\in H^s(I)$, see also
 \cite{Boyd1}. Recently, in \cite{Wang}, the author studied the speed of convergence of the expansion of such a function in a basis of PSWFs.  We should mention that the methods
used in the previous three references are heavily based on the use
of the properties of the PSWFs as eigenfunctions of the
differential operator $L_c,$ given by (\ref{eq1.0}). They pose the problem of
 the best choice of the value
of the band-width $c>0,$ for
approximating well a given $f\in H^s(I)$, but their answer is mainly experimental.  It has
been numerically checked in \cite{Boyd1, Wang} that the smaller
the value of $s,$ the larger the  value of $c$ should be.

Our study tries to give a satisfactory answer to this important problem
of the choice of the parameter $c.$
More precisely, we show that 
if  $f\in
H^s(I)$, for some positive real number $s>0,$ then for any
integer $N\geq 1,$ we have
\begin{equation}
\| f-S_N f\|_{L^2(I)}\leq K(1+c^2)^{-s/2} \|
f\|_{H^s(I)}+ K\sqrt{ \lambda_N(c)} \|f\|_{L^2(I)}.
\end{equation}
Here, ${\displaystyle S_N f=\sum_{k=0}^N <f,\psi_{n,c}> \psi_{n,c}}$ and   $K$ is a constant depending only on $s.$
Moreover, we study an  $L^2(I)-$convergence rate of the projection $S_N f$ to $f.$ This is done by using the 
 decay of the eigenvalues $(\lambda_n(c))_n$ as well as the use of some estimates of Legendre expansion coefficients of PSWFs, combined with 
the following exponential decay rate of the PSWFs expansion coefficients for  the exponential trigonometric functions
\begin{equation}
|\langle e^{ik\pi x}, \psi_{n,c}(x)\rangle|\leq M' e^{-an},\quad |k|\leq n/M,\quad n\geq \max\left( c M, 3\right).
\end{equation}
Here, $c\geq 1,$ $M\geq 1.40$ and $M', a >0$ are two positive constants. Under these hypotheses and notations, our rate of convergence
of $S_N f$ to $f\in H^s(I), s>0,\, s\not\in \frac{1}{2} +\mathbb N,$ is given as follows
\begin{equation}
\|f- S_N (f)\|_{L^2(I)}\leq M' (1+(\pi N)^2)^{-s/2}  \| f\|_{H^s}+ M' e^{-aN}  \| f\|_{L^2}.
\end{equation}

\smallskip

This work is organized as follows. In Section 2,  we list some  estimates of the PSWFs and their associated eigenvalues $\chi_n(c)$.   In Section 3, we prove a sharp exponential decay rate
of the eigenvalues $\lambda_n(c)$ associated with the integral operator $\mathcal Q_c.$
In section 4, we first  give some useful bounds of the moments of the PSWFs, then 
we  give some practical and useful estimates of the
decay of  the Legendre expansion  coefficients of the PSWFs.
  In Section 5 we first give the
quality of approximation by the PSWFs in the set  of almost time
and band-limited functions. Then, we combine these results with
those of  Section 4 and give a first  $L^2(I)-$error bound of
approximating a function $f\in H^s(I)$ by its $N$th terms
truncated  PSWFs series expansion. The proof of this bound is
based on the use of the quality of approximation of almost
bandlimited functions by the PSWFs.
Then, we study a more elaborated error analysis of the spectral approximation
by the PSWFs in the periodic Sobolev space.
This quality of approximation is then extended to the usual Sobolev space $H^s(I).$
 These new estimates provide us with a way for
the choice of the appropriate bandwidth $c>0$ to be used by a
PSWFs based method for the approximation in a given Sobolev space
$H^s(I).$ In Section 6,  we provide the reader
with  some numerical examples that illustrate the different
results	 of this work.

 We will frequently  skip the parameter  $c$ in $\chi_n(c)$ and $\psi_{n,c}$, when there is no doubt on the value of the bandwidth.
 We then note  $q= c^2/\chi_n$ and skip both parameters $n$ and $c$ when their values are obvious from the context. 

\section{Estimates of PSWFs and eigenvalues $\chi_n(c).$}
 Here we first list   some classical as well as some recent results on PSWFs and their eigenvalues $\chi_n(c)$,
 then we push forward the methods and adapt them to our study. We systematically use the same notations as in \cite{Bonami-Karoui2}.
It is well known that  the eigenvalues $\chi_n$ satisfy the classical inequalities
\begin{equation}
\label{bounds1-chi}
n (n+1)   \leq \chi_n \leq n(n+1)+c^2.
\end{equation}
In case where $q=c^2/\chi_n \leq 1,$ the following  better lower bound of $\chi_n$ has been recently given in \cite{Bonami-Karoui1},
\begin{equation}
\label{bounds2-chi}
n (n+1)+ (3-2\sqrt 2)c^2   \leq \chi_n.
\end{equation} 
Next, we recall the elliptic Legendre integral of the first and second kind,that are  given respectively, by
\begin{equation}
\label{elliptic_integrals}
\mathbf K(k)=\int_0^1 \frac{ dt}{\sqrt{(1-t^2)(1-k^2 t^2)}},\quad \quad
 \mathbf E(k)=\int_0^1 \sqrt{\frac{1- k^2 t^2}{1-t^2}}\, dt, \quad  0\leq k\leq 1.
\end{equation}
 Osipov has proved in \cite{Osipov} that the condition $q=\frac{c^2}{\chi_n}<1$
 is fulfilled when $c<\frac{\pi n}2$, while it is not when $c> \frac{\pi (n+1)}2.$ This is part of the following statement, which gives precise lower and upper bounds of the
 quantity ${\displaystyle \sqrt{q}= \frac{c}{\sqrt{\chi_n}}}$, see
 \cite{Bonami-Karoui2}.
\begin{lemma} \label{chi-between2}
For all $c>0$ and $n\geq 2$ we have
\begin{equation}\label{ineqPhi}
\Phi \left(  \frac {2c}{\pi(n+1)}\right) < \frac c{\sqrt{\chi_n}} < \Phi \left(\frac {2c}{\pi n}\right),
\end{equation}
where $\Phi$ is the  inverse of the function $k\mapsto \frac k{\mathbf E(k)}=\Psi(k),\,\, 0\leq k\leq 1.$
\end{lemma}
This is equivalent to the fact that
\begin{equation}\label{withE}
 \frac {\pi n}{2\E(\sqrt q)}<\sqrt{\chi_n}<  \frac {\pi(n+1)}{2\E(\sqrt q)}.
\end{equation}
The left hand side is due to Osipov \cite{Osipov}.
Note that $\Phi(0)=0$ and $\Phi(1)=1.$ Also, we should mention that since
\begin{equation}\label{derivé_Psi}
\Psi'(x)= \frac{\mathbf E(x) - x \mathbf E'(x)}{(\mathbf E(x))^2}=\frac{\mathbf K(x)}{(\mathbf E(x))^2},\quad 0\leq x <1,
\end{equation}
then
\begin{equation}\label{bounds_Phi'}
0\leq \Phi'(x)=\frac{(\mathbf E(\Phi(x)))^2}{\mathbf K(\Phi(x))}\leq \frac{(\mathbf E(0))^2}{\mathbf K(0)}=\frac{\pi}{2},\quad 0\leq x <1.
\end{equation}
Hence, $\Phi$ is an increasing function on $[0,1].$ Moreover, since ${\displaystyle \frac 2\pi\leq\frac{1}{\mathbf E(x)}\leq 1,}$  we have
$$\frac {2x}{\pi}\leq \Psi(x)\leq x.$$
One gets the following useful bounds of $\Phi,$
\begin{equation}\label{bounds_Phi}
x\leq \Phi(x)\leq \frac{\pi x}{2} , \quad 0\leq x\leq 1.
\end{equation}


We will use bounds for $\psi_{n, c}$ given in \cite{Bonami-Karoui2}, which have been established under the condition that $(1-q)\sqrt{\chi_n}>3.5 \E(\sqrt q)$.
Compared to  \cite{Bonami-Karoui2}, where the condition $(1-q)\sqrt{\chi_n(c)}>3.5 \E(\sqrt q)$ is systematically used to develop the uniform estimates over $[-1,1]$ of the $\psi_n$, we leave some flexibility for the choice of the constant $\kappa$. We will only need estimates at $1$, which we give here in a slightly different form compared to \cite{Bonami-Karoui2}.

Let us first recall some notations.
\begin{equation}\label{notations}
\epsn=\left((1-q)\sqrt{\chi_n}\right)^{-1}, \qquad \alpha=1. 5, \qquad \beta = 0. 37.
\end{equation}
At this moment we do not systematically replace $\alpha$ and $\beta$ by numerical values to simplify further improvements.
\begin{lemma} \label{equivalence}
Let  $n\geq 3$. We assume that the condition
\begin{equation}\label{versusChi}
(1-q)\sqrt {\chi_n(c)}>\kappa
\end{equation} is satisfied for some $\kappa\geq 4$. Then,  there exists a constant $\delta(\kappa)$ (independent of $c$ and $n$) such that  one has  the following bounds for $A=\psi_{n, c}(1)\chi_n(c)^{-1/4}$.
\begin{equation}\label{boundsA}
 \frac{\pi}{2 \mathbf  K(\sqrt{q})}\left(1-\delta(\kappa) \, \varepsilon_n\right) \leq A^2 \leq  \frac{\pi}{2 \mathbf  K(\sqrt{q})}
 \left(1+\delta(\kappa)\, \varepsilon_n\right).
\end{equation}\end{lemma}
We refer to  \cite{Bonami-Karoui1}, Theorem 3, for the  proof. Explicit values for the constant $\delta (\kappa)$  can also be deduced from \cite{Bonami-Karoui2}. We can choose
\begin{equation}\label{delta}
\delta(\kappa)=\eta\left( 2+\frac \eta{\kappa}\right), \qquad \eta = C(\kappa)\left(\frac{\beta }{1+(1- \kappa^{-1}\beta)^{1/2}}+\frac{\sqrt 2 \alpha \kappa}{\kappa-\alpha}\right)
\end{equation}
with $C(\kappa)^{-1}= (1- \kappa^{-1}\beta)^{1/2}-\frac{\sqrt 2 \alpha}{\kappa-\alpha}.$

In any case, we see that the theoretical values of $\delta(\kappa)$  are larger that $4. 6$. This corresponds to a systematic error in the approximation of the PSWFs. We find approximatively $\delta(4)\approx 90$, $\delta(12)\approx 7.7$.  Numerical tests (see Example  1 in Section 6) tend to prove that the quantities $\delta(\kappa)$  may be taken far much smaller.

\smallskip

 We have proved in \cite{Bonami-Karoui1} that one has the inequality
\begin{equation}
|A| =|\psi_{n, c}(1)|\chi_n(c)^{-1/4}\leq 2\qquad \mbox{for} \quad c\leq \frac{\pi(n+1)}{2}.
\end{equation}
So in particular the right hand side bound of \eqref{boundsA}  is not accurate when $\kappa$ is small.
\smallskip

 We need to translate Condition \eqref{versusChi} in terms of the parameters $n,c$, which can be done by using
 [Proposition 4, \cite{Bonami-Karoui2}]. The inequality given there is the following.
For $n\geq 2$ and $q<1,$
\begin{equation}\label{NtoCh}
(1-q)\sqrt{\chi_n}\geq \frac {(n-\frac{2c}{\pi})-e^{-1}}{\log n +5},
\end{equation}
A further improvement of the previous inequality is given by the following lemma:
\begin{lemma}\label{comparison} Let $n\geq 3$, $q<1$ and $\kappa\geq 4$. Then  one of the following conditions,
\begin{equation}
c\leq n-\kappa,
\end{equation}
\begin{equation}\label{NtoCh2}
\frac{\pi n} {2}-c> \frac{\kappa }4 (\ln (n)+9),
\end{equation} implies the inequality  \eqref{versusChi}, that is,
\begin{equation*}
(1-q)\sqrt {\chi_n(c)}>\kappa.
\end{equation*}
Moreover, if we assume already that $c>\frac {n+1}2$, then the condition $\frac{\pi n} {2}-c> \frac{\kappa }4 (\ln (n)+6)$ is sufficient.
\end{lemma}
\begin{proof}
Let $\gamma=\frac{2c}{\pi n}$. It follows from \eqref{withE} that
\begin{equation}\label{intermediate}
1-\gamma < 1-\sqrt q+\frac{ \E(\sqrt q)-1}{\E(\sqrt q)}.
\end{equation}
We claim that
\begin{equation}
\label{EEq2}
\mathbf E(x)-1 \leq  (1-x^2)\left( \frac{1}{4}\ln \left(\frac 1{1-x^2}\right)+\ln 2\right).
\end{equation}
Let us assume this and go on with the proof. It follows that
\begin{equation}
1-\gamma< \frac{1-q}{\E(\sqrt q)} \left( \frac{1}{4}\ln \left(\frac 1{1-q}\right)+\frac {\E(\sqrt q)}{1+\sqrt q}+\ln 2\right).
\end{equation}
We then use the elementary inequality, valid for $0<s<1$,
$$ s\ln ( 1/s)\leq 1/n + s\ln (n/e).$$
It implies that $$1-\gamma-\frac 1{4n \E(\sqrt q)}<\frac{1-q}{\E(\sqrt q)} \left( \frac{1}{4}\ln (n/e)+\frac {\E(\sqrt q)}{1+\sqrt q}+\ln 2\right).$$
We use also   \eqref{withE} to conclude that
\begin{equation}\label{versus1-q}
(1-q)\sqrt{\chi_n}\geq \frac{\pi n}{2\E(\sqrt q)}(1-q)>\kappa
\end{equation}
whenever
$$\frac{\pi n} {2}-c> \kappa \left(\frac 14 \ln (n/e)+)\frac {\E(\sqrt q)}{1+\sqrt q}+\ln 2\right)+ \frac 1{4n}.$$
This is the case, in particular, when $\frac{\pi n} {2}-c> \frac \kappa 4 \left(\ln (n)+9\right)$, using the fact that $\frac {\E(\sqrt q)}{1+\sqrt q}\leq \frac \pi 2$.

The condition $c\geq \frac {n+1}2$ implies that $q>\frac 1\pi$. Then, by  using the value of $\E(\sqrt {\pi^{-1}})$, the constant $9$ in (\ref{NtoCh2}) can be replaced by $6$.

It remains to prove \eqref{EEq2}. We write
\begin{eqnarray}
\mathbf E(x)-1 &\leq & (1-x^2)\int_0^1\frac 1{(\sqrt{1-x^2t^2}+\sqrt{1-t^2})}\frac{t\, dt}{\sqrt{1-t^2}}\label{intermediate2}\\
&=&\int_0^1\frac {ds}{(1-x^2+s^2x^2)^{\frac 12}+s}.
\end{eqnarray}
We cut the last integral into two parts. For the first one, from $\sqrt{1-x^2}$ to $1$, we replace the denominator by $2s$ and find the logarithm term. For the second one we replace the denominator by $\sqrt{1-x^2}+s$ and find $\ln 2$.
\end{proof}

We will need another inequality of the same type:
\begin{equation}\label{comparisonK}
1-\frac{2c}{\pi n}\leq 2(1-q)\K(\sqrt q).
\end{equation}
This is a consequence of \eqref{intermediate}, using the fact that $\E(x)-1\leq (1-x^2)\K(x)$, which comes directly from \eqref{intermediate2}.

\medskip
We end this section by giving  bounds for the values of the successive derivatives of $\psi_n$ at  $0$. We have proved in \cite{Bonami-Karoui1} that
\begin{equation}
\label{bound1psi}
 |\psi_n(0)|^2+ \chi_n^{-1}|\psi'_n(0)|^2\leq \left\{\begin{array}{ll}
 1  &\mbox{ if\ \  } 0\leq q \leq 2\\ \frac{q+1}{\sqrt q} &\mbox{ if\ \  } q>  2.\end{array}\right.
\end{equation}
Let us prove that, furthermore, successive derivatives at $0$ may be also bounded.

\begin{proposition} Assume that ${\displaystyle q=\frac{c^2}{\chi_n}<1}$.  Then for any  integer $ k\geq 0$ satisfying $k(k+1)\leq \chi_n,$ we have
\begin{equation}\label{ineqderivatives}
\left|\psi^{(k)}_{n}(0)\right|\leq (\sqrt{\chi_n})^k.
\end{equation}
\end{proposition}
\begin{proof} Because of \eqref{bound1psi}, it is sufficient to prove that $m_k=(\sqrt{\chi_n})^{-k}\left|\psi^{(k)}_{n}(0)\right|$ is bounded by $(|\psi_n(0)|^2+ \chi_n^{-1}|\psi'_n(0)|^2)^{1/2}$. Moreover, since $\psi_{n, c}$ has same parity as $n$, then it is sufficient to consider even derivatives or odd derivatives depending on the parity of $n$. Assume first that  $n$ is even and consider $k=2l$.
  We show that for a fixed $n,$ ${\displaystyle
\psi^{(2l)}_{n,c}(0)}$ has alternating signs, that is
$\psi^{(k)}_{n,c}(0)\psi^{(k-2)}_{n,c}(0)<0.$ Indeed, by an
iterative use  of the identity
$$(1-x^2)\psi_n''(x)=2x \psi_n'(x)+(c^2x^2-\chi_n)\psi_n(x),$$
 one can easily check that the
${\displaystyle \psi^{(k)}_{n,c}(0)=\psi^{(k)}(0)}$ are given by
the following recurrence relation,
\begin{equation}\label{recurrence1}
\psi^{(k+2)}(0)=(k(k+1)-\chi_n) \psi^{(k)}(0)+k(k-1) c^2
\psi^{(k-2)}(0),\quad k\geq 0,
\end{equation}
with $\psi(0) >0,$ $\psi^{(2)}(0)=-\chi_n \psi(0).$ Note that
$\psi^{(2)}(0)\psi(0)<0.$ Assume that
$\psi^{(k)}(0)\psi^{(k-2)}(0)<0.$ Multiplying both sides of (\ref{recurrence1}) by $\psi^{(k)}(0),$ using the assumption that
$k(k+1)\leq \chi_n$ as well as the induction hypothesis, one
concludes that the induction assumption holds for the order $k.$
Consequently, we have,
\begin{equation}\label{recurrence2}
\left|\psi^{(k+2)}(0)\right|=(\chi_n-k(k+1))
\left|\psi^{(k)}(0)\right|+k(k-1) c^2
\left|\psi^{(k-2)}(0)\right|,\quad k\geq 0.
\end{equation}
This may be rewritten as
\begin{equation}\label{recurrence3}
m_{k+2}=\left(1-\frac{k(k+1)}{\chi_n}\right)m_k+k(k-1)\frac{q}{\chi_n}
m_{k-2}.
\end{equation}
The fact that all $m_{2l}$ are bounded by $m_0=|\psi(0)|= m_2$ follows at once by induction.
For $n$ odd the proof follows the same lines. \end{proof}

As a consequence of the previous proposition, we have the following corollary concerning the sign and the bounds
of the different moments of the $\psi_n.$

\begin{corollary} Let $c>0,$ be a positive real number. We assume that $q= c^2/\chi_n<1$. Then, for $j(j+1)\leq \chi_n$, all moments $\int_{-1}^1 y^j \psi_{n}(y)\, dy$ of the same parity as $n$ have the same sign and
\begin{equation}
\label{momentspsi}
\left|\int_{-1}^1 y^j \psi_{n}(y)\, dy\right| \leq \left(\frac{1}{q}\right)^{j/2}|\mu_n(c)|.
\end{equation}
\end{corollary}

\begin{proof}
By taking the $j-$th derivative at zero on both sides of ${\displaystyle \int_{-1}^1 e^{icxy} \psi_{n}(y)\, dy =\mu_n(c) \psi_{n}(x),}$
one gets
\begin{equation}
\label{moments2psi}
\int_{-1}^1 y^j \psi_{n}(y)\, dy= (-i)^j c^{-j} \mu_n(c)\psi_{n}^{(j)}(0),\quad\mbox{with } i^2=-1.
\end{equation}
Since $\psi_{n}^{(j)}(0)$ and  $\psi_{n}^{(j+2)}(0)$ have opposite signs, then the previous equation implies
that  moments  have the same sign for any positive integer $j$ with $j(j+1)\leq \chi_n.$ The second inequality
of (\ref{momentspsi}) follows from the previous proposition.
\end{proof}

\section{Sharp decay estimates of eigenvalues $\lambda_n(c).$}

In this section, we use some of the estimates we have given in the previous section and
we prove a sharp  over-exponential decay rate of  the eigenvalues $(\lambda_n(c))_n.$
We first recall that these $\lambda_n(c)$ are governed by the following differential equation,
see for example \cite {Xiao},
\begin{equation}\label{eq1lambda_n}
\partial_c \ln \lambda_n(c)=\frac{2|\psi_{n,c}(1)|^2}{c}.
\end{equation}
As a consequence, for fixed $n$ there exists a unique value of $c$ for which $\lambda_n(c)=1/2$, which we call $c_n^*$. We know  from \cite{Landau} that it can be bounded below and above, namely
 \begin{equation}
 \label{validityC0}
 \frac{\pi}{2}(n-1)\leq c^*_n\leq \frac{\pi}{2} (n+1)\quad\mbox{ with}\quad \lambda_n(c^*_n)=\frac{1}{2}.
 \end{equation}
 By combining (\ref{eq1lambda_n}) and (\ref{validityC0}), one gets
 \begin{equation}
 \label{eqq2lambda_n}
 \lambda_n(c)=\frac{1}{2} \exp\left(-2\int_c^{c^*_n} \frac{(\psi_{n,\tau}(1))^2}{\tau} \, d\tau \right).
  \end{equation}

Our main result is the following theorem.
\begin{theorem}\label{Required0}
There exist three constants $\delta_1\geq 1, \delta_2, \delta_3, \geq 0$ such that, for $n\geq 3$ and $c\leq \frac{\pi n}2$,
\begin{equation}\label{required0}
\delta_1 ^{-1}n^{-\delta_2}\left(\frac c{c+1}\right)^{\delta_3}\leq \frac{\widetilde { \lambda_n(c)}}{ \lambda_n(c)}\leq \delta_1 n^{\delta_2}\left(\frac c{c+1}\right)^{-\delta_3},
\end{equation}
where  \begin{equation}\label{decay2}
\widetilde { \lambda_n(c)} = \frac{1}{2} \exp\left(-\frac{\pi^2(n+\frac 12)}{2} \int_{\Phi\left(\frac{2c}{\pi (n+\frac 12)}\right)}^1 \frac{1}{t(\mathbf E(t))^2}\, dt \right).\end{equation}
\end{theorem}
The factor $\frac c{c+1}$ can be replaced by $1$ when $c>1$ and replaced by $c$ when $c<1$. We have written the formula this way  to avoid to have to distinguish between the two cases, $c\geq 1$ and $0<c<1.$

  It is simpler to write equivalent  inequalities for logarithms, which is done in the following proposition. We keep the same notations for constants, which are of course not the same. We note $\ln^+ (x)$  the positive part of the Logarithm, that is, $\max(0, \ln (x))$. The following theorem
  is a fundamental theorem that is required in the proof of the main Theorem \ref{Required0}.
\begin{theorem}\label{Required}
There exist three non negative constants $\delta_1, \delta_2, \delta_3$ such that, for $n\geq 3$ and $c\leq \frac{\pi n}2$, we have
\begin{equation}\label{required}
\int_c^{c^*_n} \frac{(\psi_{n,\tau}(1))^2}{\tau} \, d\tau =\frac{\pi^2(n+\frac 12)}{4} \int_{\Phi\left(\frac{2c}{\pi (n+\frac 12)}\right)}^1 \frac{1}{t(\mathbf E(t))^2}\, dt  +\mathcal  E,
\end{equation}
with
\begin{equation}
|\mathcal E|\leq \delta_1+\delta_2 \ln (n)+\delta_3 \ln^+ (1/ c).
\end{equation}\end{theorem}

Let us make some comments before starting the proof. At this moment the three constants are not sufficiently small and cannot be used reasonably to obtain  numerical values. But they can be computed and are not that enormous. There is no hope, of course, to have found an exact formula for $\lambda_n(c)$ and \eqref{decay2} gives only an approximation. But these theoretical approximation errors may be seen as a kind of theoretical validation of the quality of approximation, which we test numerically in Section 6.
\smallskip

It has been observed by many authors, and predicted by the work of Landau and Widom \cite{LW}, that for fixed $c$ the eigenvalues $ \lambda_n(c)$ decrease first exponentially in some interval starting at $ [\frac{2c}{\pi}]+1$ with length a multiple of $\ln(c)$, then super-exponentially as in the asymptotic behavior given by Widom. This is what one observes in Formula \eqref{decay2}, but the error terms do not allow to observe the  decay rate
at the small decay starting region. In fact the tools that we use, that is, the lower and upper bounds for $\psi_{n,\tau}(1)^2$, are only valid for $c_n^* -\tau$ sufficiently large in terms of $\ln(n)$.

We try to have  small constants at each step but are certainly far from the best possible. We give an explicit bound for $\mathcal E$ in
\eqref{final}.

The following notations will be used frequently in the sequel. We define
\begin{eqnarray}
I(a, b)&=&\int_a^{b} \frac{(\psi_{n,\tau}(1))^2}{\tau} \, d\tau.\\
\mathcal J (y) &=&\frac{\pi^2 }{4} \int_{\Phi\left(\frac{2y}{\pi }\right)}^{1} \frac{1}{t(\mathbf E(t))^2}\, dt
\end{eqnarray}

We should mention that the proofs of Theorems 2 and Theorem 3, require many steps, so
 we start by giving a sketch of these proofs.

\medskip

\noindent{\sl \bf Sketch of the proofs.}
\smallskip

We want to prove that
$$I(c, c_n^*)\approx (n+\frac 12)\mathcal J \left(\frac c{n+\frac{1}{2}}\right).$$
Lemma \ref{equivalence} expresses the fact that, under some condition depending on a parameter $\kappa$, we have
$$\psi_{n,\tau}(1))^2\approx \frac{\pi\sqrt{\chi_n(\tau)}}{2 \mathbf  K(\sqrt{q(\tau)})}=\frac{\pi\tau}{2 \sqrt{q(\tau)}\mathbf  K(\sqrt{q(\tau)})}.$$
The parameter $\kappa$ is related with the quality of approximation and Lemma \ref{comparison} proves that the condition for this   may be written $c<c_n^\kappa$ for some $c_n^\kappa$. From the last equivalence, it follows that
$$I(c, c_n^\kappa)\approx \int_c^{c_n^\kappa} \frac{\pi d\tau}{2 \sqrt{q(\tau)}\mathbf  K(\sqrt{q(\tau)})} .$$
Then Lemma \ref{chi-between2} will be interpreted as the fact that
$$\sqrt{q(\tau)}\mathbf  K(\sqrt{q(\tau)})\approx \Phi\left(\frac {2\tau}{\pi(n+\frac 12)}\right)\mathbf  K\circ\Phi\left(\frac {2\tau}{\pi(n+\frac 12)}\right).$$
It is then elementary to rely the new integral with the function $\mathcal J$ and finally find that
$$I(c, c_n^\kappa)\approx (n+\frac 12)\mathcal J \left(\frac c{n+\frac{1}{2}}\right).$$ It remains to bound the tails of the integrals  $I( c_n^\kappa, c_n^*)$, which we can do because the two values are sufficiently close.
\smallskip

Let us start the proof itself. We need a set of intermediate results that can be classified into three main steps.
The first step will concern the properties of the function $\mathcal J$. In the second step,  we give bounds of the  tails of the integrals.
Finally, in the third step, we use the results of the previous two steps and complete the proofs of  Theorems 2 and 3.
\medskip

\noindent{\sl \bf First step: Properties of $\mathcal J$.}
\smallskip

We define \begin{equation}\label{Jl}
\mathcal J_l(c)=\frac \pi 2\int_c^{\frac{\pi l}2}\frac{d\tau}{\Phi\left(\frac {2\tau}{\pi l}\right)\mathbf  K\circ\Phi\left(\frac {2\tau}{\pi l}\right)}.
\end{equation}
Such integrals are clearly involved in the proof as seen in the sketch. We first see that they are related with $\mathcal J$.
\begin{lemma}
We have the identity
\begin{equation}
\mathcal J_l(c)= l\mathcal J(c/l). \label{identity}
\end{equation}
\end{lemma}
\begin{proof}
We  consider   the substitution
\begin{equation}
\label{subs1}
 s = \Phi\left(\frac{2\tau}{\pi l}\right), \qquad   \tau= \frac{\pi l}2 \Psi(s).
\end{equation}
We have already seen in \eqref{derivé_Psi} that
$ \Psi'(x)=\frac{\mathbf K(x)}{(\mathbf E(x))^2}.$
Hence, we have
$$ \mathcal J_l(c)= l\int_{\Phi(\frac{2c}{\pi l})}^1\frac{ ds }{s (\mathbf E(s))^2}= l \mathcal J(c/l).$$
\end{proof}

The following proposition gives us  upper and lower  bounds, as well as the asymptotic behavior of $\mathcal J$.
 \begin{proposition}
For $x\in (0, \pi/2)$, one has the upper and lower bounds
\begin{equation}\label{majJ}
\ln^+\left(\frac{1}{ x}\right)\leq \mathcal J(x)\leq \frac{\pi^2 }{4} \ln\left(\frac{\pi}{2 x}\right).
\end{equation}
Moreover, one can write
   \begin{equation}\label{at-infty}
 \mathcal J(x)=\frac{\pi^2 }{4}\int_{\Phi(2x/\pi)}^1\frac{dt}{t(\mathbf E(t))^2}= \ln\left(\frac{4}{ex}\right)+ \mathcal E',
 \end{equation}
 with $|\mathcal E'|\leq \frac{\pi^2 x^2}{8} $.
\end{proposition}
\begin{proof}
The first inequalities are an easy consequence of bounds below and above of $\Phi,$ given by
(\ref{bounds_Phi}).
Let us prove \eqref{at-infty}. We first write, for $0<y<1$,
\begin{equation}
\frac{\pi^2 }{4}\int_y^1 \frac{dt}{t(\mathbf E(t))^2} +\ln (y)=\Delta-\int_{0}^y \frac{\frac{\pi^2}{4}-\E(t)^2}{t(\mathbf E(t))^2} \,
 dt =\Delta - I_1(y). \end{equation}
Here
$$\Delta= \int_{0}^1 \frac{\frac{\pi^2}{4}-\E(t)^2}{t(\mathbf E(t))^2} \, dt.$$
It is probably well-known that
\begin{equation}\label{Delta}
\Delta=\ln\left(\frac{4}{e}\right)
\end{equation}
but we did not find any reference. We will see it as a corollary of Widom's Theorem.
The integral $I_1(y)$ is bounded by $\frac{\pi^2 y^2}8$. This  is a consequence of the elementary inequalities
$$1\leq \mathbf E(s)\leq \frac{\pi}{2},\quad \frac \pi 2 -\E(s)\leq x^2\int_0^1\frac{t^2 \,dt}{\sqrt{1-t^2}}= \frac{\pi s^2} 4.$$
Let us now fix $y=\Phi(2x/\pi)$. At this point we have proved that
$$0\leq \ln \left( \frac x y \right)-\mathcal E'= I_1(y) \leq \frac{\pi^2 y^2}8.$$
From the inequalities
$$\frac{2y}\pi\leq\frac{2x}\pi=\frac y{\E(y)}\leq \frac{2y}\pi (1-\frac{y^2}{2})^{-1}\leq \frac{2y}\pi (1+y^2),$$
it follows that $0\leq  \ln \left( \frac x y \right)+y^2$. We have proved the proposition.
\end{proof}
 This proposition leads to the following corollary, where we recognize the equivalent given by Widom.
\begin{corollary}\label{Tilde}
We have the double inequality
\begin{equation}\label{tilde}
\frac{1}{2}\left(\frac{ec}{4(n+\frac 12)}\right)^{2n+1}e^{-\frac {\pi^2}4 \frac{c^2}{n+\frac 12}}\leq \widetilde { \lambda_n(c)}\leq \frac{1}{2}\left(\frac{ec}{4(n+\frac 12)}\right)^{2n+1}e^{+\frac {\pi^2}4 \frac{c^2}{n+\frac 12}}.
\end{equation}
\end{corollary}

\begin{proof}
Just note that ${\displaystyle \widetilde { \lambda_n(c)}=\frac{1}{2}\exp\left(-(2n+1)\mathcal J(c/(n+1/2))\right) }$ and use (\ref{at-infty})
with $x=\frac{c}{n+1/2}.$
\end{proof}
Let us go back to quantities $\mathcal J_l$. It is a straightforward consequence of \eqref{identity}  that the quantity $\mathcal J_l(c)$ increases with $l$. The next lemma gives reverse inequalities.
\begin{lemma}
 We have the inequalities
 \begin{equation}\label{comparisonJ}
\mathcal J_{n+1}(c)-\frac {\pi^2}{8} \ln \left(\frac{\pi(n+1)}{2c}\right)-\frac {\pi^3}{16} \leq \mathcal J_{n+\frac 12}(c)\leq \mathcal J_n(c)-\frac {\pi^2}{8} \ln \left(\frac{\pi(n+\frac 12)}{2c}\right)+\frac {\pi^3}{16}.
 \end{equation}
\end{lemma}
\begin{proof}
We will prove only one of the inequalities, the other one being identical. Elementary computations give
$$ \mathcal J_{n+1}(c)-\mathcal J_{n+\frac 12}(c)\leq \frac 12 \mathcal{J}\left(\frac{c}{n+1}\right)+\frac {\pi^2}4 (n+\frac 12)
\ln\left(\frac{\Phi\left(\frac{2 c}{\pi (n+\frac 12)}\right)}{\Phi\left(\frac{2 c}{\pi (n+1)}\right)}\right).$$
We use \eqref{majJ} for the first term.  The second one is bounded by
$$
\frac {\pi^2}4 (n+\frac 12)\frac{\Phi\left(\frac{2 c}{\pi (n+\frac 12) }\right)-\Phi\left(\frac{2 c}{\pi (n+1)}\right)}{\Phi\left(\frac{2 c}{\pi (n+1)}\right)}\leq \frac{\pi^3}{16}.$$
Indeed, this is a consequence of the fact that $\Phi'(x)\leq \pi/2$ and ${\displaystyle \frac{x}{\Phi(x)}\leq 1,\,}$ for $ 0<x\leq 1$.
\end{proof}

\medskip

\noindent{\bf Second step: tails of the integrals.}
\smallskip

We fix some constant $\kappa\geq 4$ (for instance $\kappa=12$) and we assume that $n\geq 2\kappa+1$. Then, we know from Lemma \ref{comparison},  that the condition \eqref{versusChi}, that is,
$$(1-q)\sqrt {\chi_n}>\kappa,$$ is satisfied for $c<\frac{n+1}{2}$.
Next, if we  define
 \begin{equation}
\label{validityC}
  c^\kappa_n =\max\left(\frac{\pi n}{2}  -\frac \kappa 4 (\ln (n)+6), \frac{n+1}{2}\right)
\end{equation}
then, we have the following.
\begin{lemma}
For $n\geq 2\kappa +1$ we have the inequality
\begin{equation}\label{Ilarge}
I(c_n^\kappa, c_n^*) \leq \pi \kappa \ln (n)+ 6\pi\kappa+2\pi^2.
\end{equation}
\end{lemma}
\begin{proof}
Recall that ${\displaystyle |\psi_{n,c}(1)|\leq 2 \chi_n^{1/4}}$ and  ${\displaystyle \sqrt{\chi_n(c)}\leq \frac{\pi}{2}(n+1)}$,
\cite{Osipov},   so that
$$|\psi_{n, \tau}(1)|^2\leq 4\sqrt{\chi_n(\tau)}\leq 2\pi(n+1).$$
Hence, we have
\begin{eqnarray*}
\int_{c^\kappa_n}^{c^*_n} \frac{(\psi_{n,\tau}(1))^2}{\tau} \, d\tau &\leq& 2\pi(n+1) \ln \left (1+\frac {\frac \pi 2+\frac\kappa 4(\ln (n)+6)} {c^\kappa_n}\right).
\end{eqnarray*}
We conclude by using the fact that $c^\kappa_n\geq \frac{n+1}{2}$.
 \end{proof}
We claim that we can conclude the proof of Theorem \ref{Required} when $n\geq 2\kappa +1$ and $c <c_n^{\kappa}$. More precisely, we get, under these conditions, the inequalities
\begin{equation}
-\frac{\pi^2}{16}( \kappa \ln (n)+ 6\pi\kappa+\pi)\leq I(c, c^*_n)-(n+\frac 12)\mathcal J\left(\frac{c}{n+\frac 12}\right)\leq \pi \kappa \ln (n)+ 6\pi\kappa+2\pi^2
\end{equation}
The right hand side comes from the previous lemma, the left hand side from \eqref{majJ}.

We will conclude this paragraph by showing that we have also the conclusions of Theorem \ref{Required} and Theorem \ref{Required0} for the finite number of missing values of $n$, that is, $n\leq 2\kappa +1$. There is no problem to have upper bounds and lower bounds that do not depend on $c$ for $c<1$. From Corollary  \ref{Tilde}, we have a precise estimate in terms of $c^{2n+1}$ for $\widetilde{\lambda_n(c)}$. The same is given for $\lambda_n(c)$ by the following lemma.
\begin{lemma}
 Assume that $n\geq 1$ is fixed and let $0<c<1$. Then there exist two constants $\delta(n), \delta'(n)$ such that
 \begin{equation}\label{chi}
 \delta(n) \, c^{2n+1}\leq \lambda_n(c)\leq  \delta'(n)\,  c^{2n+1}.
 \end{equation}
\end{lemma}
\begin{proof}
We first note that  $I(1, c^*_n)\leq I(1, \frac{\pi(n+1)}2)$. We recall that on this interval we have the inequality $|\psi_{n, \tau}(1)|^2\leq 4\frac{\pi(n+1)}2$. So $I(1, c^*_n)\leq 2\pi (n+1) \ln(\frac{\pi(n+1)}2)$.
Inside the integral defining $I(c, 1)$
we use  the following inequality, that may be found in \cite{Bonami-Karoui2},
\begin{equation}
\left| |\psi_{n, \tau}(1)|-\sqrt{n+\frac 12}\right|\leq \frac{\tau^2}2.
\end{equation}
So $\left|I(c, 1)-(n+\frac 12)\ln\left(\frac 1c\right)\right|\leq 1$, from which we conclude.
\end{proof}

\smallskip

\noindent{\bf Third step: Proofs of Theorems 2 and 3.}
\smallskip

We fix $\kappa>4$.  Because of the previous steps, which allowed to conclude in the other cases, we can assume that
$$ n\geq 2\kappa +1 \qquad \qquad c<c_n^{\kappa}= \max\left(\frac{\pi n}{2}  -\frac \kappa 4 (\ln (n)+6), \frac{n+1}{2}\right).$$
In view of \eqref{required}, we want to give a bound to
$$\mathcal E=I(c, c_n^*)-\left(n+\frac 12\right)\mathcal J\left(\frac{c}{n+\frac 12}\right).$$
We have already given a bound to a first error term
$$\mathcal E_1=I(c, c_n^*)-I(c, c_n^\kappa).$$
Because of \eqref{Ilarge} we know that
\begin{equation}\label{E1}
0\leq \mathcal E_1\leq \pi \kappa \ln (n)+ 6\pi\kappa+2\pi^2.
\end{equation}
Next, the conditions on $\kappa$ allow to use the double inequality \eqref{boundsA}. Namely,
\begin{equation}
 \label{psi1}
 \left(\psi_{n,\tau}(1)\right)^2 = \frac{\pi}{2 \mathbf K(\sqrt{q})} \sqrt{\chi_n(\tau)} +\mathcal{R}(\tau),\quad
 |\mathcal{R}(\tau)|\leq \frac{\delta(\kappa )}{(1-q(\tau)) \mathbf K(\sqrt{q(\tau)})},\quad 0\leq \tau \leq c^{\kappa}_n.
 \end{equation}
This leads to a second error, 
$$\mathcal E_2=I(c,c_n^{\kappa})-\frac{\pi}{2}\int_c^{c_n^{\kappa}}\frac{d\tau}{\sqrt{q(\tau)} \mathbf K(\sqrt{q(\tau)})},$$ which is bounded by
$$|\mathcal E_2|\leq\delta(\kappa)\int_c^{c_n^\kappa}\frac 1{(1-q(\tau)) \mathbf K(\sqrt{q(\tau)})}\frac{ d\tau}\tau.$$
 \begin{lemma}
 We have the inequality
 \begin{equation}
 |\mathcal E_2|\leq 2\delta (\kappa)\left((1+\frac{\pi \kappa}{4} ) \ln(n)+\ln^+(c)+\frac{3\pi \kappa}{2}\right).
 \end{equation}
 \end{lemma}
 \begin{proof}
By \eqref{comparisonK}, we know  that
 $$2(1-q(\tau)) \mathbf K(\sqrt{q(\tau)})\geq 1- \frac{2 \tau}{\pi n}.$$
 So we have the inequality
 $$ |\mathcal E_2|\leq 2\delta (\kappa) \int_{\frac{2 c}{\pi n}}^{\frac{2 c_n^\kappa}{\pi n}}\frac{ ds}{(1-s)s}
 \leq 2\delta (\kappa)\left(\ln \left(\frac nc\right)+ \ln\left(\frac 1{1-\frac{2 c_n^\kappa}{\pi n}}\right)\right).$$
We conclude at once.
\end{proof}

It remains to consider the main term, that is,
\begin{equation}
I_{\rm main}(c, c_n^\kappa)=\frac{\pi}{2}\int_c^{c_n^\kappa}\frac{\sqrt{\chi_n(\tau)}}{ \mathbf K(\sqrt{q(\tau)})} \frac{d\tau}{\tau}=\frac{\pi}{2}\int_c^{c_n^\kappa}\frac{d\tau}{ \sqrt {q(\tau)}\mathbf K(\sqrt{q(\tau)})} .
\end{equation}
We will use the monotonicity properties of $\sqrt {q(\tau)}\mathbf K(\sqrt{q(\tau)})$ , namely
$$\Phi\left(\frac{2 \tau}{\pi (n+1)}\right)\mathbf K\circ\Phi\left(\frac{2 \tau}{\pi (n+1)}\right)\leq\sqrt {q(\tau)}\mathbf K(\sqrt{q(\tau)})
 \leq \Phi\left(\frac{2 \tau}{\pi n}\right)\mathbf K\circ\Phi\left(\frac{2 \tau}{\pi n}\right).$$
It follows that
$$ \mathcal{J}_{n}(c)- \mathcal{J}_{n}(c_n^\kappa)\leq I_{\rm main}(c, c_n^\kappa)\leq \mathcal{J}_{n+1}(c).$$
So the last error, 
$$\mathcal E_3=I_{\rm main}(c, c_n^\kappa)- \mathcal{J}_{n+\frac 12}(c)=I_{\rm main}(c, c_n^\kappa)- \left(n+\frac 12\right)\mathcal{J}\left(\frac{c}{n+\frac{1}{2}}\right),$$ satisfies the inequalities
$$\mathcal{J}_{n}(c)-\mathcal{J}_{n+\frac 12}(c)- \mathcal{J}_{n}(c_n^\kappa)\leq\mathcal E_3\leq \mathcal{J}_{n+1}(c)-\mathcal{J}_{n+\frac 12}(c).$$
It remains to use \eqref{comparisonJ} and \eqref{majJ} to conclude. We finally find that
\begin{equation}\label{final}
| \mathcal E | \leq \pi \kappa \ln (n)+ 6\pi\kappa+2\pi^2+2\delta (\kappa)\left((1+\frac{\pi \kappa}{4} ) \ln(n)+\ln^+(c)+\frac{3\pi \kappa}{2}\right)+\frac {\pi^2}{8} \ln \left(\frac{\pi(n+\frac 12)}{2c}\right)+\frac {\pi^3}{16}.
 \end{equation}
 So we can take the following values for $\delta_1, \delta_2, \delta_3$.
 \begin{eqnarray*}
 \delta_1&=&22+3\pi\kappa(2+\delta(\kappa))\\
 \delta_2&=&\frac {\pi^2}{8}+\pi \kappa+2\delta (\kappa)(1+\frac{\pi \kappa}{4} )\\
 \delta_3&=& \frac {\pi^2}{8}+2\delta (\kappa)(1+\frac{\pi \kappa}{4} ).
 \end{eqnarray*}

It is easy to see from the proof above that this bound is also valid for $c_n^\kappa<c<\frac{n\pi}2$, that is, under the assumptions of
Theorem \ref{Required}, except for the values $n\leq 2\kappa+1$. These estimates are not sharp enough to justify a further study to minimize the sum by a specific choice of $\kappa$. When $\kappa=12$ we find $\delta_2\approx 200$. We could have improved bounds at each step, but not significantly. Numerical experiments tend to prove that they are much smaller.
\smallskip

This concludes the proofs of Theorem \ref{Required} and Theorem \ref{Required0}.

\bigskip

 From Theorem \ref{Required} and Corollary \ref{Tilde} we get the following corollary:
 \begin{corollary} There exist three constants $\delta_1\geq 1, \delta_2, \delta_3, \geq 0$ such that, for $n\geq 3$ and $c\leq \frac{\pi n}2$,
 \begin{equation}\label{best}
A(n, c)^{-1}\left(\frac{ec}{2(2n+1)}\right)^{2n+1}\leq  \lambda_n(c)\leq A(n,c)\left(\frac{ec}{2(2n+1)}\right)^{2n+1}.
\end{equation}
 with $$A(n, c)=\delta_1 n^{\delta_2}\left(\frac c{c+1}\right)^{-\delta_3}e^{+\frac {\pi^2}4 \frac{c^2}{n}}.$$
 \end{corollary}

Widom's Theorem says that $A(n, c)$ can be replaced by a quantity that tends to $1$ for $n$ tending to $\infty$. We cannot give such an asymptotic behavior at this moment, but we can estimate errors for fixed $c$ and $n$, which he does not. Remark that we used the fact that $\Delta= \ln(4/e)$, see (\ref{Delta}), without proving it or giving a reference. This is a consequence of the asymptotic behavior found by Widom, which cannot  be valid at the same time as \eqref{best} if $e/4$ is replaced by another constant. 
 This implies in particular Theorem \ref{th-intro}.
 
 \smallskip
 
 It may be useful to give also the following  corollary.
 \begin{corollary} There exist   constants $a>0$ and $\delta\geq 1$ such that, for $c\geq 1$ and $n>1.35\, c$, we have  
 \begin{equation}\label{best-bis}
 \lambda_n(c)\leq \delta e^{-an}.
\end{equation}
\end{corollary}
\begin{proof} The constant $1.35$ has been chosen so that $2\ln (\frac{4 n}{ec})>\frac {\pi^2c^2}{4n^2}$ for $n>1.35 \, c$.
\end{proof}

Also, by using (\ref{best}), one gets the following corollary.
 \begin{corollary} Let $c\geq 1,$ then for any $0\leq a <\frac{4}{e},$ there exists $N_a\in \mathbb N$ such that 
 $$\lambda_n(c)\leq e^{-2n\log\left(\frac{a n}{c}\right)},\quad\forall\,\, n\geq N_a.$$ 
 Moreover, for any $b>\frac{4}{e},$ there exists $N_b\in \mathbb N$ such that 
  $$\lambda_n(c)>  e^{-2n\log\left(\frac{b n}{c}\right)},\quad\forall\,\, n\geq N_b.$$ 
\end{corollary}

\section{Decay estimates of the Legendre expansion coefficients}

Recently, there is an extensive amount of work devoted to
new highly accurate computational methods of the PSWFs, see
\cite{Beylkin,Karoui1,  Xiao}. In particular, the methods
given in \cite{Beylkin, Xiao} are based on an efficient quadrature
method on the unit circle that provides highly accurate values of
the PSWFs  inside $[-1,1],$ as well as accurate approximations of
the different eigenvalues $\mu_n(c),\, n\geq 0.$ The methods
developed in \cite{Karoui1} for computing the values of
the $\psi_{n,c}(x)$ inside $[-1,1]$ and the eigenvalues $\mu_n(c)$
are  based on an appropriate matrix representation of the
finite Fourier transform operator $\mathcal F_c,$ given by (\ref{eq1.1}). Also, we should mention
a classical method known as Flammer's method, \cite{Flammer} that uses the differential operator $\mathcal L_c$, is extensively used
to compute the PSWFs and their eigenvalues.
This method  is based on the following  Legendre expansion of the PSWFs,
\begin{equation}\label{eqqq3.2}
 \psi_{n}(x)={\sum_{k\geq 0}} \beta_k^n \overline{P_k}(x).
 \end{equation}
Recall that  $\psi_n$ has the same parity as $n.$ Hence, the previous Legendre expansion coefficients of the 
$\psi_n$ satisfy $\beta_k^n=0$ if $n$ and $k$ have different parities.
The expansion (\ref{eqqq3.2}) is in particular used to compute the eigenvalues in terms of the coefficients $\beta_k^n$. Indeed, using the fact that the Fourier transform of the Legendre polynomials can be expressed  in terms of Bessel functions, as well as the property of $\psi_n$ of being an eigenfunction of $\mathcal{F}_c$, we have
\begin{equation}\label{eq2.2.10}
\psi_{n}(x)= \frac{\sqrt{2\pi}}{|\mu_n(c)|}{\sum_{k\geq 0}} (-1)^k
\beta_k^n \sqrt{k+1/2}\frac{J_{k+1/2}(c
x)}{\sqrt{cx}},
\end{equation}
which extends analytically outside $I$. Here $J_{\alpha}$ denotes the Bessel function of the first kind and order $\alpha > -1$.
As a consequence, Slepian has proved in \cite{Slepian1} that
\begin{equation}\label{eq2.2.11}
 \mu_n(c)=\frac{2\pi }{c} \left[\frac{{\sum_{k\geq 0}} i^k \sqrt{k+1/2}
\,\,\beta_k^n \,\, J_{k+1/2}(c)}{ {\sum_{k\geq 0}}
\beta_k^n\sqrt{k+1/2}}\right],
\end{equation}
is the exact value of the $n-$th eigenvalue of the finite Fourier
transform operator $\mathcal F_c.$

 It is well known that the  different
expansion coefficients $(\beta_k^n)_k$ as well as the
corresponding eigenvalues $\chi_n$ are obtained by solving the
following eigensystem
\begin{eqnarray}
&&\frac{(k+1)(k+2)}{(2k+3)\sqrt{(2k+5)(2k+1)}} c^2 \beta_{k+2}^n
 + \big( k(k+1)+ \frac{2k(k+1)-1}{(2k+3)(2k-1)} c^2\big)
\beta_k^n \label{eigensystem}\\
&&\hspace*{3cm} + \frac{k(k-1)}{(2k-1)\sqrt{(2k+1)(2k-3)}} c^2
\beta_{k-2}^n= \chi_n(c) \beta_k^n, \quad k\geq 0.\nonumber
\end{eqnarray}

A  useful decay estimate of the $\beta_k^n$ is based on the following
 positivity result of the  $\beta_k^n.$

 \begin{lemma}
Let $c>0,$ be a fixed positive real number. Then, for all positive  integers $k, n$ such that $k(k-1)+1. 13\, c^2\leq \chi_n(c)$, we have $\beta_k^n\geq 0$.
\end{lemma}

\begin{proof}
We recall that the $\beta_k^n$ are given by the  eigensystem (\ref{eigensystem}). Let us first consider $k=2$ (when $n$ is even) and $k=3$ (when $n$ is odd) satisfying the condition. We compute
$$ \beta_2^{n} =\frac{3\sqrt{5}}{2 c^2}\left(\chi_n -\frac{c^2}{3}\right)\beta_0^{n}\geq \beta_0^{n}\geq 0,\quad
 \beta_3^{n} =\frac{5\sqrt{21}}{6 c^2}\left(\chi_n-2 -\frac{3 c^2}{5}\right)\beta_1^{n}\geq \beta_1^{n}\geq 0.$$
  For $k\geq 2$, taking upper bounds for the fractions as in  \cite{Chen}, Equation (\ref{eigensystem}) implies that
$$ \frac{2c^2}{3\sqrt{5}} ( \beta_{k+2}^n +\beta_{k-2}^n)\geq 
 ( \chi_n(c)-k(k+1)-\frac{11c^2}{21}) \beta_k^n.
$$
The constant $1. 13$ has been chosen so that $1.13>\frac{4}{3\sqrt{5}}+\frac{11}{21}$. By the assumption on $k$  it is not possible for $\beta_{k+2}$ and $\beta_{k-2}$ to be bounded by $\beta_k$, with $\beta_k>0$. The same is valid at each  step $k-2, k-4,\cdots$. Since $\beta_2\geq \beta_0 >0$ (resp.  $\beta_3\geq \beta_1 >0$) depending on the parity of $n$,  the sequence $\beta_{2j}$ (resp. $\beta_{2j+1})$ is non decreasing for $2j\leq k+2$ (resp. $2j+1\leq k+2$). This implies the positivity.
\end{proof}

The following proposition provides us with a useful decay rate of the expansion coefficients $\beta_k^n.$

\begin{proposition}
Let $c>0,$ be a fixed positive real number. Then, for all positive integers  $n,k$ such that $k(k-1)+1. 13\,  c^2\leq \chi_n(c)$, we have
\begin{equation}
\label{Decay2beta}
|\beta_0^n|\leq \frac{1}{\sqrt{2}}|\mu_n(c)|\quad \mbox{ and }\quad |\beta_k^n|\leq \sqrt{\frac{5}{4\pi}}\left(\frac{2}{\sqrt{q}}\right)^{k} |\mu_n(c)|.
\end{equation}
\end{proposition}

\begin{proof}  The first inequality follows from Corollary 1 and   the fact that
${\displaystyle \beta_0^n=\frac{1}{\sqrt{2}}\int_{-1}^1 \psi_n(y)\, dy.}$ To prove the second inequality, we first note that the moments of the normalized Legendre polynomials are non-negative and they are given in \cite{Andrews}, by
\begin{equation}\label{eq2.3.2}
a_{jk}=\int_{-1}^{1}x^{j}\overline{P_k}(x)dx=\left\{\begin{array}{ll}
0 &\mbox {\ if\  } j< k \mbox{ or }
 j-k \mbox{ is odd}\\
\frac{\sqrt{\pi}j!}{2^{j}\left(\frac{j-k}{2}\right)!\Gamma(\frac{k+j+3}{2})}&\mbox
{\ if\  } j-k \geq 0 \mbox{ and  } j-k \mbox{ is even.}
\end{array}\right.
\end{equation}
Since ${\displaystyle x^j=\sum_{k=0}^j a_{jk} \overline{P_k}(x),}$ then the moments of the $\psi_n$ are related to the
PSWFs Legendre expansion coefficients by the following rule,
$$\int_{-1}^1 x^j \psi_n(x)\, dx = \sum_{k=0}^j a_{jk}  \beta_k^n.$$
Since by the previous lemma, we have $\beta_k^n \geq 0,$ for any $0\leq  k\leq j$ and since the $a_{jk}$ are positive, then the previous equality
implies that
\begin{equation}\label{decay3beta}
\beta_j^n \leq \frac{1}{a_{jj}} \int_{-1}^1 x^j \psi_n(x)\, dx \leq  \frac{1}{a_{jj}} \left(\frac{1}{q}\right)^{j/2}|\mu_n(c)|.
\end{equation}
The last inequality follows from the previous corollary. On the other hand, we have
$$a_{jj}= \frac{\sqrt{\pi} \sqrt{j+1/2} j! }{2^j \Gamma(j+3/2)}= \frac{\sqrt{\pi}  j! }{2^j \sqrt{j+1/2}\Gamma(j+1/2)}.$$
Moreover, it is well known that ${\displaystyle j^{1-s}\leq \frac{\Gamma(j+1)}{\Gamma(j+s)}\leq (j+1)^{1-s}}.$ Hence, we have
\begin{equation}\label{decay4beta}
\frac{1}{a_{jj}}\leq \frac{2^j}{\sqrt{\pi}}\sqrt{1+\frac{1}{2j}}\leq 2^j\sqrt{\frac{5}{4\pi}},\quad\forall\, j\geq 1.
\end{equation}
By combining (\ref{decay3beta}) and (\ref{decay4beta}), one gets the second inequality of (\ref{Decay2beta}).
\end{proof}

\begin{remark}\label{condition2prop3}  The condition $k(k-1)+1.13 c^2\leq \chi_n(c)$ of the previous proposition can be replaced with the following  more explicit condition.
Consider a real number  $A>1,$ then  by using  (\ref{bounds2-chi}), one concludes that  if $n\geq c \alpha \frac{ A}{\sqrt{A^2-1}}$ with $\alpha=\sqrt{1.13+2\sqrt{2}-3}\approx 0.979$ and $k\leq n/A,$ then  the conditions for \eqref{Decay2beta}  are satisfied.
Moreover, 
the previous constant $0.979$  is certainly not optimal, it is a  consequence of the non optimal  lower 
bound of $\chi_n(c),$ given by (\ref{bounds2-chi}). 
\end{remark}

\begin{remark} In \cite{Chen}, by using the eigensystem  (\ref{eigensystem}), the authors have  obtained under a stronger condition,
a decay of the Legendre coefficients $(\beta_k^n)_{k}$ which is similar but less precise to 
 the one we have given by (\ref{Decay2beta}). More precisely, they have shown that if $0< k\leq 2m$ with 
$m=O(n^{2/3})$ and $2m(2m+1)<\frac{\log 2}{2} \chi_n,$ then there exists a constant $D$ such that 
$$|\beta_k^n| \leq D\left(\frac{2}{\sqrt{q}}\right)^k |\beta_0^n|,\,\,\mbox{ for even }k,\quad   |\beta_k^n| \leq D\left(\frac{2}{\sqrt{q}}\right)^k |\beta_1^n|, \mbox{ for odd } k.$$
\end{remark}

\section{Quality of the spectral approximation by the PSWFs}

In this section, we first study the quality of approximation of
almost band-limited functions by the classical PSWFs, $\psi_{n}$
that are  concentrated  on $[-b,b],$ for some $b>0.$ Then, we
extend this study to the case of periodic and non periodic Sobolev
space $H^s([-1,1]), s>0.$

\subsection{Approximation of almost time and band-limited functions}

In this paragraph, $\|\cdot\|_2$ denotes the norm in $L^2(\R)$.
We show that the set $\{\psi_{n}(x),\,\,
n\geq 0\}$ is well adapted for the representation of almost
time-limited and almost band-limited functions, which are defined
as follows.

\begin{definition}
 Let $T=[-a,+a]$ and $\Omega=[-b, +b]$ be two intervals. A
function  $f$, which we assume to be normalized in such a way that
$\|f\|_2=1$, is said to be $\epsilon_T-$concentrated in $T$ and
$\epsilon_{\Omega}-$band concentrated in $\Omega$ if
$$\int_{T^c} |f(t)|^2\, dt \leq \epsilon_T^2,\qquad
\frac 1{2\pi}\int_{\Omega^c} |\widehat f(\omega)|^2\, d\omega \leq
\epsilon_{\Omega}^2.$$
\end{definition}

Up to a re-scaling of the function $f$, we can always assume that
$T=[-1, 1]$ and $\Omega=[-c, +c]$, with $c:=ab$. Indeed, for $f$
that is $\epsilon_T-$concentrated in $T=[-a, +a]$ and
$\epsilon_{\Omega}-$band concentrated in $\Omega=[-b, +b]$, the
normalized function $g(t)=\sqrt a f(at)$ is
$\epsilon_T-$concentrated in $[-1, +1]$ and
$\epsilon_{\Omega}-$band concentrated in $[-ab, +ab]$.

Before stating the theorem, let us give some notations. For $f\in L^2(\mathbf R),$
 we consider its expansion $f=\sum_{n\geq
0} a_n \psi_{n,c}$ in $L^2([-1,+1]).$ Due to the
normalization of the functions $\psi_{n,c}$ given by \eqref
{eqq1.4}, the following equality holds,
\begin{equation}\label{plancherel-1+1}
   \int_{-1}^{+1}|f(t)|^2 dt=\sum_{n\geq 0}|a_n|^2.
\end{equation}
 We call $S_{N,c}f,$  the $N$-th partial sum, defined by
\begin{equation}\label{plancherel}
   S_{N,c}f(t)=\sum_{n<N}a_n \psi_{n,c}(t).
\end{equation}
 We write more simply $S_Nf$ when there is no ambiguity. In the next lemma, we prove that $S_N f$ tends to $f$ rapidly
 when $f$ belongs to the space of band-limited functions. This statement is both very simple and classical, see for instance \cite{Rokhlin2, Slepian1}
 or Theorem 3.1 in  \cite{Wang}.

 \begin{lemma}\label{band-limited} Let $f\in B_c$ be an $L^2$ normalized function. Then
 \begin{equation}\label{eq6.0}
 \int_{-1}^{+1}|f- S_Nf|^2 dt\leq   \lambda_N(c).
\end{equation}
\end{lemma}
 \begin{proof} Since the set of functions $\psi_{n,c}$ is also an orthogonal basis of $B_c$, the function $f$ may be written on $\R$ as $f=\sum_{n\geq 0} a_n \psi_{n,c}$, with
  \begin{equation}\label{plancherelR}
   \int_{\R}|f(t)|^2 dt=\sum_{n\geq 0}|\lambda_n(c)|^{-1}|a_n|^2.
\end{equation}
The two expansions coincide on $[-1, +1]$, and, from
\eqref{plancherelR} applied to $f-S_N f$, it follows that
$$\int_{-1}^{+1}|f- S_Nf|^2 dt\leq  \sup_{n\geq N} |\lambda_n(c)| \sum_{n\geq N}|\lambda_n(c)|^{-1}|a_n|^2.$$
We use the fact that the sequence $ |\lambda_n(c)|$ decreases and
\eqref{plancherelR} to conclude.
\end{proof}

Next we define  the time-limiting operator $P_T$ and the
band-limiting operator $\Pi_{\Omega}$ by:
$$ P_T(f) (x) = \chi_T(x) f(x),\qquad
\Pi_{\Omega}(f)(x)=\frac{1}{2\pi}\int_{\Omega} e^{i x
\omega}\widehat{f}(\omega)\, d\omega.$$ The following proposition
provides us with the quality of approximation of almost time- and
band-limited functions by the PSWFs.
\begin{proposition}
If $f$ is an $L^2$ normalized function that is
$\epsilon_T-$concentrated in $T=[-1, +1]$ and
$\epsilon_{\Omega}-$band concentrated in $\Omega=[-c, +c],$ then
for any positive integer $N,$ we have
\begin{equation}\label{eqqq6.0}
\left( \int_{-1}^{+1}|f- S_Nf|^2 dt\right)^{1/2} \leq  \epsilon_{\Omega}+\sqrt{\lambda_N(c)}
 \end{equation}
and, as a consequence,
\begin{equation}\label{eeqq6.0}
\|f - P_T S_N f\|_2 \leq
\epsilon_T+\epsilon_{\Omega}+\sqrt{\lambda_N(c)}.
\end{equation}
More generally, if $f$ is an $L^2$ normalized function that is
$\epsilon_T-$concentrated in $T=[-a, +a]$ and
$\epsilon_{\Omega}-$band concentrated in $\Omega=[-b, +b]$ then,
for $c=ab$ and for any positive integer $N,$ we have
\begin{equation}\label{eeqq6.01}
\|f - P_T S_{N,c,a}f\|_2 \leq \epsilon_T+\epsilon_{\Omega}+\sqrt
{\lambda_N(c)}
\end{equation}
where $S_{N,c,a}$ gives the $N$-th partial sum for the orthonormal
basis $\frac 1{\sqrt a}\psi_{n,c}(t/a)$ on $[-a, +a]$.
\end{proposition}

\begin{proof}[{\bf Proof:}] We first  prove \eqref{eqqq6.0} by writing $f$ as the sum of $\Pi_\Omega f$ and $g$.
Remark first that $\int_{-1}^{+1}|g-S_N g|^2 dt\leq \|g\|_2 \leq \epsilon_\Omega$. We then use Lemma \ref{band-limited}
for the band limited function $\Pi_\Omega f$ to conclude. The rest of the proof follows at once.
\end{proof}

\begin{remark} Let $f$ be a normalized $L^2$ function that vanishes outside $I$ and  we assume that
 $f\in H^s(\R)$.
Then $f$ gives an example of $0$-concentrated in $I$ and $\epsilon_c$-band concentrated in $[-c, +c]$,
with $\epsilon_c\leq M_f/c^s$ and ${\displaystyle M_f^2=\frac 1 {2\pi}\int |\widehat f(\xi)|^2|\xi|^{2s} d\xi.}$
\end{remark}

\subsection{Approximation by the PSWFs in Sobolev spaces}

In this paragraph, we study the quality of approximation by the
PSWFs in the Sobolev space $H^s([-1,1]).$  We
 provide an $L^2([-1,1])$-error bound of the approximation
of a function  $f\in H^s([-1,1])$ by the $N-$th partial sum of its
expansion in the basis of PSWFs.

To simplify notation we will write $I=[-1,1]$. We should mention that
different spectral approximation results by the PSWFs in $H^s(I)$ have
been already given in \cite{Boyd1, Chen, Wang}.  More precisely,
the following result has been proved in \cite{Chen}. Here $a_k
(f)=\int_{-1}^1 f(x)\psi_{k}(x)\, dx$.

\begin{theorem}\mbox{$({\bf Theorem\, 3.1\, in }\,\cite{Chen}).$} Let $f\in H^s(I),\,  s\geq 0$. Then
$$|a_N(f)|\leq C \left( N^{-2/3 s}\| f\|_{H^s(I)}+\left(\sqrt{\frac{c^2}{\chi_N(c)}}\right)^{\delta N}\| f\|_{L^2(I)}\right),$$
where $C, \delta$ are independent of $f, N$ and $c.$
\end{theorem}



In \cite{Wang}, the author has used a different approach for the
study of the spectral approximation by the PSWFs. More precisely,
by considering the weighted Sobolev space $\widetilde H^r(I),$
associated with the differential operator ${\mathcal L}_c$
defined by
$$\widetilde H^r(I)= \left\{ f\in L^2(I),\, \|f\|^2_{\widetilde H^r(I)}=\|{\mathcal L}_c^{r/2} f\|^2=\sum_{k\geq 0}
(\chi_k)^r | f_k|^2 <+\infty\right\},$$ where $f=\sum f_k$ is the expansion in the basis of PSWFs.  The following
result has been given in \cite{Wang}.

\begin{theorem}\mbox{$({\bf Theorem\, 3.3\, in } \, \cite{Wang})$}.
For any $f\in \widetilde H^r(I),$ with $r\geq 0,$ we have
$$\|f-S_N f\|_{L^2(I)} \leq (\chi_N(c))^{-r/2} \|f\|_{\widetilde H^r(I)}\leq N^{-r} \|f\|_{\widetilde H^r(I)}.$$
\end{theorem} 

It is important to mention that the  error bounds of the spectral
approximations given by the previous two theorems, do not indicate
how to choose a    ``good" value of the bandwidth $c$ 
to approximate a given $f\in H^s(I).$ By a simultaneous use of the
properties of the PSWFs as eigenfunctions of the differential
operator $\mathcal L_c$ and the integral operator $\mathcal F_c,$ we give a first
answer to this question. This is the subject of the following
theorem.

\begin{theorem}
Let  $c > 0$ be a positive real number. Assume that $f\in
H^s(I)$, for some positive real number $s>0$. Then for any
integer $N\geq 1,$ we have
\begin{equation}\label{eq222.1}
\| f-S_N f\|_{L^2(I)}\leq K(1+c^2)^{-s/2} \|
f\|_{H^s(I)}+ K\sqrt{ \lambda_N(c)} \|f\|_{L^2(I)}.
\end{equation}
Here, the constant $K$  depends only on $s$. Moreover it can be
taken equal to $1$ when $f$ belongs to the space $H^s_0(I)$.
\end{theorem}

\begin {proof} To prove  (\ref{eq222.1}), we first use the fact that for
any real number $s\geq 0,$ there exists a linear and continuous
extension operator $E: H^s(I)\rightarrow H^s(\mathbb R).$
Moreover, if $f\in H^s(I)$ and $F= E(f) \in H^s(\mathbb R),$
then there exists a constant $K>0$ such that
\begin{equation}\label{extension}
\| F \|_{L^2(\mathbb R)} \leq K \| f\|_{L^2(I)},\qquad \|
F\|_{H^s(\mathbb R)} \leq K \|f\|_{H^s(I)}.
\end{equation} We recall that the Sobolev norm of a function F on $\R$ is given by
$$\|F\|_{H^s(\R)}^2 =\frac{1}{2\pi}\int_{\R} (1+|\xi|^2)^s|\widehat f (\xi)|^2\, d\xi.
$$
In particular, for $F$ $c-$bandlimited, one has
$$\|F\|_{L^2(\R)}^2\leq (1+c^2)^{-s}\|F\|_{H^s(\R)}^2.$$
Next,
if ${\mathcal F}$ denotes the Fourier transform operator and if
$${\cal G} = {\cal F}^{-1}(\widehat F \cdot 1_{[-c,c]}),\quad {\cal
H}={\cal F}^{-1}(\widehat F \cdot (1-1_{[-c,c]})),$$ then ${\cal
G}$ is $c-$bandlimited and $F={\cal G} + {\cal H}.$ Moreover,
since $\| \widehat {\mathcal G}\|_{L^2(\mathbb R)}\leq \|\widehat
F\|_{L^2(\mathbb R)}$ and $\|  {\mathcal H}\|_{L^2(\mathbb R)}\leq
c^{-s} \|F\|_{H^s(\mathbb R)},$ then by using (\ref{extension}),
one gets
\begin{equation}\label{eq222.3}
\| {\mathcal G}\|_{L^2(\mathbb R)}\leq K \|f\|_{L^2(I)},\quad \|
{\mathcal H}\|_{L^2(I)}\leq K (1+c^2)^{-s/2} \| f\|_{H^s(I)}.
\end{equation}
Finally, by using the previous inequalities and the fact that
${\mathcal G}$ is $c-$bandlimited, one concludes that
\begin{eqnarray*}
\| f- S_N f\|_{L^2(I)}&\leq & \|\mathcal G - S_N \mathcal G\|_{L^2(I)}
+\|\mathcal H - S_N \mathcal H\|_{L^2(I)}\\
&\leq & \sqrt{ \lambda_N(c)} \|{\mathcal
G}\|_{L^2(\mathbb R)} +\|\mathcal H\|_{L^2(I)}\\
&\leq &  \sqrt{ \lambda_N(c)} K \|f\|_{L^2(I)} +
 K (1+c^2)^{-s}\| f\|_{H^s(I)}.
\end{eqnarray*}

This concludes the proof for general $f$. When $f$ is in the subspace $H^s_0(I)$, one can take as
extension operator the extension by $0$ outside $I$, so that the constant $K$ can be replaced by $1$.
\end{proof}

\begin{remark}
This should be compared with the results of \cite{Wang}, given by
Theorem 5. This has the advantage to give an error term for all
values of $c$, while the first term in \eqref{eq222.1} is only
small for $c$ large enough. On the other hand, Wang compares his
specific Sobolev space with the classical one and finds that
$$\|f\|_{\widetilde{H^s}(I)}\leq C(1+c^2)^{s/2}\|f\|_{H^s(I)}.$$
For large values of $N$ we clearly have  ${\displaystyle \frac{(1+c^2)}{\chi_N}\ll
(1+c^2)^{-1},}$ but it goes the other way around when $\chi_N$ and
$1+c^2$ are comparable. So it may be useful to have both kinds of
estimates in mind for numerical purpose and for the choice of the
value of  $c.$
\end{remark}

\begin{remark}
 The error bound given by the previous theorem has the advantage to be explicitly given in terms of
 $c$ and $\lambda_n(c).$ Nonetheless, it has a drawback that it does not imply a rate of convergence, nor the  convergence
 of $S_N(f)$ to $f$ in the usual $L^2(I)-$norm. To overcome this problem, we devote the remaining of this section
 to a more elaborated convergence analysis in the $2$-periodic Sobolev space $H^s_{per},$ then we extend this analysis
 to the usual $H^s(I)-$space.
\end{remark}

Next, we  consider  the subspace $H^s_{per}$ of functions in $H^s(I)$ that extend into $2-$periodic functions of the same regularity.
For such functions, one can also use the norm
$${\displaystyle \|f\|_{H^s_{per}}=\sum_{k\in \mathbb Z}
(1+(k\pi)^2)^s | b_k(f)  )|^2.}$$
 Here, $$ b_k(f)=\frac 1{\sqrt 2} \int_{-1}^{+1}f(x) e^{-i\pi k x} dx=\frac 1{\sqrt 2}\widehat f (k\pi) $$
is the coefficient of the Fourier series expansion of $f.$ We then have the following theorem.

 \begin{theorem}
Let $c\geq 1,$   then there exist constants $M>1.40$ and $ M',\, a>0$ such that, when  $N\geq \max (c M,3)$ and  $f\in H^s_{per}, s>0$,  we have the inequality
\begin{equation}\label{error1}
\|f- S_N (f)\|_{L^2(I)}\leq M' (1+(\pi N)^2)^{-s/2}  \| f\|_{H^s_{per}}+ M' e^{-aN}  \| f\|_{L^2}.
\end{equation}
\end{theorem}
\begin{proof}  We start with reductions of the problem, which are analogous to the ones that we have detailed above. It is sufficient to prove this separately with the constant $M'/2$ for periodic functions $g$ and $h=f-g$, where $g$ is the projection of $f$ onto the subspace of $H^s_{per}$ whose Fourier coefficients $b_k(f)$ are zero for $|k|>N/M.$ Moreover, we have directly the inequality without a second term, since the $L^2$ norm of $h$ may be  bounded by the first term multiplied by some constant. So, let us prove the inequality for $g$. This time we will prove that the inequality holds without the first  term, that is,
$$\|g- S_N (g)\|_{L^2(I)} \leq   \frac{M'}{2} e^{-aN}  \| g\|_{L^2(I)}.$$
The next reduction consists of restricting to exponentials $e^{ik\pi x}$, with $|k|\leq N/M$. Indeed, assume that we prove the previous  inequality  for all of them, with a uniform bound by $M''e^{-a'N}$. Then, by linearity we  will have
$$\|g- S_N (g)\|_{L^2(I)} \leq M''e^{-a'N}\sum | b_k(g)  )|\leq  M''e^{-a'N} \sqrt{2[N/M]+1}\;e^{-aN} \| g\|_{L^2(I)}.$$
This in turn gives constants  the required form by choosing $a<a'$. 

So we  content ourselves to consider $f(x)=e^{ik\pi x}$, with $|k|\leq N/M$. Finally, since $\|f-S_N f\|_{L^2(I)}^2=\sum \langle f, \psi_n\rangle^2$, it is sufficient to have such an estimate for each $n>N$, and conclude by taking the sum $\sum_{n> N} e^{-an}$. So the proof is a consequence of the following lemma.
\end{proof}

\begin{remark}
The previous theorem gives the rate of convergence of the truncated PSWFs series expansion of a function $f$ from  $H^s_{per}.$ This rate of convergence 
will be generalized in the sequel to the usual $H^s(I)-$space. Note that this rate of convergence drastically improves the one given by \cite{Chen}.
Moreover, unlike the error bound given in \cite{Wang}, the decay of the error bound given by the previous  theorem is still even when $N$ is comparable to $c.$ 
 Nonetheless, in practice,  Theorem 6 is useful in the sense that provides us
with a criteria for the choice of the bandwidth $c>0,$ that depends on magnitude of the Sobolev exponent $s>0.$ The smaller $s,$ the larger $c$ should be and vice versa. 
\end{remark}

\begin{lemma}
Let $c\geq 1,$ then there exist constants $M>1.40$ and $  M', \, a>0$ such that, when  $n\geq \max\left( c M, 3\right)$ and  $f(x)=e^{ik\pi x}$ with $|k|\leq n/M,$ we have
\begin{equation}
|\langle f, \psi_n\rangle|\leq M'e^{-an}.
\end{equation}
\end{lemma}

\begin{proof} This scalar product can be written  by using \eqref{eq2.2.10}
\begin{eqnarray*}
< e^{ik\pi x},\psi_n >&=& \int_{-1}^1  e^{ik\pi x}\,\psi_n(x)\, dx = \sum_{m\geq 0} \beta_m^n < e^{ik\pi x},\overline P_m > =\sum_{m\geq 0} \beta_m^n \sqrt{\frac{2}{k}}\sqrt{m+1/2} J_{m+1/2} (k\pi)\\
&=&\sum_{m= 0}^{[n/M]} \beta_m^n \sqrt{\frac{2}{k}}\sqrt{m+1/2} J_{m+1/2} (k\pi)+\sum_{m\geq [n/M]+1} \beta_m^n \sqrt{\frac{2}{k}}\sqrt{m+1/2} J_{m+1/2} (k\pi)\\
&=&I_1^n+ I_2^n.
\end{eqnarray*}
To bound $I_1^n,$ we first remark   that  the Fourier transform of $\overline P_n\chi_{[-1, 1]}$ is bounded by $1$ and then
we use remark \ref{condition2prop3} to check that \eqref{Decay2beta} is satisfied whenever $n\geq c M$ with $M\geq 1.40.$
Hence, we have 
\begin{eqnarray*}
|I_1^n|&\leq& \sum_{m= 0}^{[n/M]} |\beta_m^n| \leq  \sqrt{\frac{5}{4\pi}}|\mu_n(c)| \sum_{m= 0}^{[n/M]}
\left(\frac{2\sqrt{\chi_n}}{c}\right)^m\\
&\leq & K  \left(\frac{2\sqrt{\chi_n}}{c}\right)^{[n/M]+1} |\mu_n(c)|.
\end{eqnarray*}
Moreover, taking into account the decay of the $\mu_n(c)$ given by \eqref{best-bis} and using the upper bound
of $\chi_n,$ we conclude that
\begin{equation}
\label{boundI1}
|I_1^n| \leq K' \left(\frac{\pi(n+1)}{c}\right)^{\frac nM +1}e^{-\delta n}\leq K'' e^{-an}
\end{equation}
for some   sufficiently small positive real number $a$, as soon as $M>1.40$.
To bound $I_2^n,$ it suffices to use the fact that $|\beta_k^n|\leq 1$ and the bound
 of the Bessel function given by \cite{Andrews},
 \begin{equation}\label{Eq2.2}
 |J_{\alpha}(x)|\leq \frac{|x|^{\alpha}}{2^{\alpha}
 \Gamma(\alpha+1)},\quad \forall\, \alpha > -1/2,\quad\forall\,
 x\in \mathbb{R},
 \end{equation}
 one concludes that
 \begin{eqnarray*}
 |I_2^n|&\leq &\sum_{m\geq n/M}  \sqrt{2/k} \sqrt{m+1/2} |J_{m+1/2} (k\pi)|\leq \sum_{m\geq [n/2]+1} \sqrt{2/k}\sqrt{m+1/2}
 \frac{ (k\pi)^{m+1/2}}{2^{m+1/2} \Gamma(m+3/2)} \\
 &\leq & \sum_{m\geq [n/M]+1}\frac{ (k\pi)^m}{2^m\sqrt{m+1/2}\Gamma(m+1/2)}.
 \end{eqnarray*}
Moreover, since ${\displaystyle \Gamma(m+1/2)\geq m!/\sqrt{m+1}}$ and ${\displaystyle m! \geq (m/e)^m \sqrt{2\pi m}},$ each term is bounded by an exponential $e^{-an}$ and we find the required estimate for $|I_2^n|$.
\end{proof}

\begin{remark}
We also have a bound of the error for ordinary polynomials.
Indeed, if we consider the polynomial $f(x):=x^j$, then
$$a_n(f)=\int_{-1}^1 y^j \psi_{n,c}(y)\, dy= (-i)^j c^{-j} \mu_n(c)\psi_{n,c}^{(j)}(0),\quad\mbox{with } i^2=-1.$$
We can then use Proposition 5  to conclude that if $c^2/\chi_N
<1,$ then
\begin{equation}\label{majoration1}
\| f- S_N f\|_2^2 \leq C^2 \sum_{k\geq N} \left(\frac{\chi_k(c)}{c^2}\right)^j|\mu_k(c)|^2.
\end{equation}
\end{remark}

As a corollary of the previous theorem and remark, we obtain the following corollary that extends the result of the
previous theorem to the case of the usual Sobolev space $H^s([-1,1]).$

\begin{corollary}
Let $c\geq 1,$  and let  $s>0$ with  $[s]= m\in \mathbb N,$ and $s\not\in \frac{1}{2}+\mathbb N.$ Let  $f\in H^s(I),$ then there exist constants $M\geq 1.40$ and $M', M'_s>0$ such that, when  $ N\geq \max\left( c M, 3\right)$,   we have the inequality
\begin{equation}\label{error2}
\|f- S_N (f)\|_{L^2(I)}\leq M'_s (1+ N^2)^{-s/2}  \| f\|_{H^s([-1,1])}+ M' e^{-aN}  \| f\|_{L^2([-1,1])}.
\end{equation}
\end{corollary}

\begin{proof}
 Since $f\in H^s([-1,1])$ with $[s]= m,$ and $s\not\in \frac{1}{2}+\mathbf N,$ then there exists
a polynomial $P$, of degree  at most $m$,  such that $f+ P\in H^s_{per}.$ 
Consequently, by using the previous theorem and the inequality (\ref{majoration1}), one concludes
for (\ref{error2}).
\end{proof}

\section{Numerical results}

In this section, we illustrate the results of the previous  sections
by various numerical examples. \\

\noindent
{\bf Example 1:} In this first example, we illustrate the fact that the actual values of the  constants $\kappa$ and $\delta(\kappa),$
given  by (\ref{versusChi}) and (\ref{boundsA}), respectively, are far  much smaller than the theoretical values given in the proof of Lemma \ref{equivalence}. We are interested in these values for $n\geq 2c/\pi$.
For this purpose, we have considered the values of $c= m \pi, m=10, 20, 30,  40.$ Then, we have used Flammer's method and computed
high accurate values of $\chi_n(c)$ and $\psi_{n,c}(1).$ Then, we have  computed the smallest value of $\kappa,$ denoted by $\kappa_c$ and
ensuring the bounds  (\ref{boundsA}). Also, we have computed the corresponding values  $\delta(\kappa_c)$ so that
$A^2$ is equal to its upper bound given in (\ref{boundsA}).  It turns out that $\kappa_c,$ the critical value of $\kappa,$  is obtained for $n-$th eigenvalues $\chi_n(c)$ with
$n=n_c = [2c/\pi].$ Also, by considering various consecutive values of $n_c\leq n\leq  n_c+40$ and by computing the corresponding values of $\kappa$ and $\delta(\kappa),$ we found that the $\max \delta(\kappa)$ is of the same size as $\kappa_c.$
 Table 1 shows the  values of the critical values $\kappa_c$ and $\delta(\kappa_c)$ for the  different  values of the
bandwidth $c.$ Also, we give the values of $\max \delta(\kappa).$

\begin{center}
\begin{table}[h]
\vskip 0.2cm\hspace*{4cm}
\begin{tabular}{ccccc} \hline
 $c$ &$n_c$&$ \kappa_c$&$\delta(\kappa_c)$&$\max \delta(\kappa).$ \\   \hline
$10\, \pi$  &20    &  0.447   & 0.058  & 0.091    \\
$20\, \pi$ &40    & 0.413    & 0.051  & 0.084\\
$30\, \pi$ &60    & 0.394    & 0.047   & 0.080    \\
$ 40\, \pi$ &80   & 0.335     &0.025   & 0.048   \\
\hline
\end{tabular}
\caption{Critical values  of $\kappa,$  $\delta(\kappa)$ and $\max \delta(\kappa)$ for different values of $c.$}
\end{table}
\end{center}

\noindent
{\bf Example 2:} In this example, we compare the  explicit formula 
  given by Theorem 2 to compute highly accurate values of $\lambda_n(c)$. For this purpose, we have considered the values of $c=10 \pi, 20\pi, 30\pi$ and
computed   $\lambda_n(c)$ by using the method given in \cite{Karoui1}. Then, we have implemented
our formula (\ref{decay2}) in a Maple computing software code. Figure 1 (a), (b), (c) show the graph of  $\ln(\lambda_n(c))$ versus the graph of  $\ln(\widetilde \lambda_n(c)),$ for the different values
of $c$ and $n.$ Also, we have plotted in Figure 2, the graphs of the corresponding values of
${\displaystyle \ln\left(\frac{\lambda_n(c)}{\widetilde{ \lambda_n(c)}}\right)}.$ These figures illustrate the surprising precision of the explicit formula of Theorem 3 for computing the $\lambda_n(c)$
which is numerically valid whenever $q<1.$

\begin{figure}[h]
{\includegraphics[width=15cm,height=6.2cm]{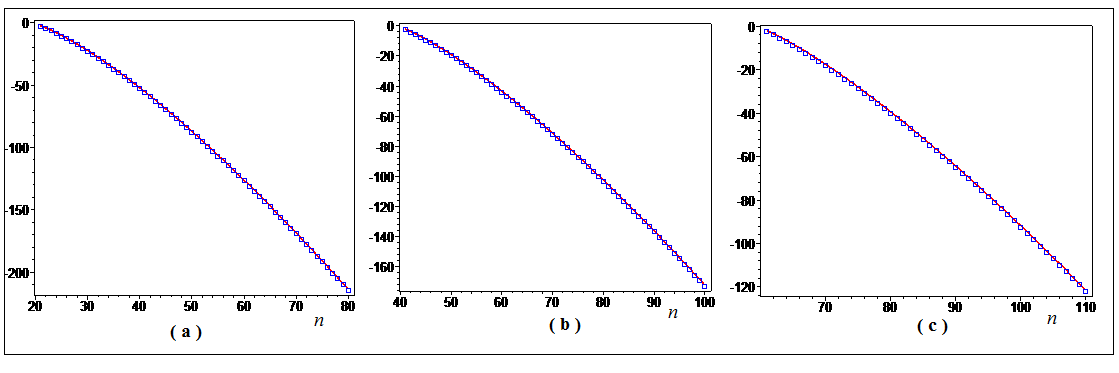}} \vskip
-0.5cm\hspace*{2cm}\caption{Graphs of
$\ln (\widetilde{ \lambda_n(c)})$ (boxes) and $\ln (\lambda_n(c))$ (red) with $c=10 \pi$ for (a), $c=20\pi$ for  (b) and $c=30\pi$ for (c). }
\end{figure}

\begin{figure}[h]
{\includegraphics[width=15cm,height=5cm]{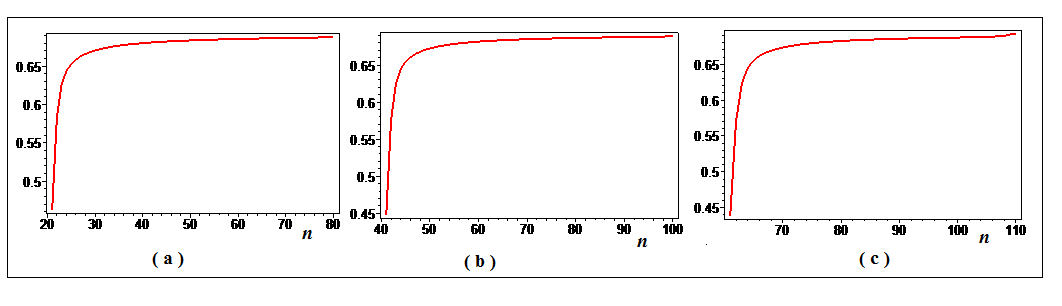}} \vskip
-0.5cm\hspace*{2cm}\caption{ Graphs of
$\ln\left(\frac{\lambda_n(c)}{\widetilde{ \lambda_n(c)}}\right)$ with $c=10\pi$ for (a), $c=20\pi$ for (b) and $c=30\pi$ for (c).}
\end{figure}

\vskip 0.5cm
\noindent
Next, to illustrate the quality of approximation by the $\psi_n$
in the Sobolev space $H^s(I),$ we first describe a numerical method for the
computation of the  PSWFs series
expansion coefficients of a function from the Sobolev space
$H^s(I).$ Note that if $f\in H^s_{per},\, s>0,$ then its different
PSWFs series expansion coefficients can be easily approximated as follows.
For a positive
integer $K,$  an approximation $a_n^K(f)$ to $a_n(f)$ is given by
the following formula
\begin{equation}\label{eq3.10}
a_n^K(f)= \frac{\mu_n(c)}{\sqrt{2}} \sum_{k=-K}^K b_k(f)
\psi_{n,c}\left(\frac{k\pi}{c}\right)= a_n(f) + \epsilon_K,
\end{equation}
where the $b_k(f)$ are the  Fourier coefficients of $f$ and
where ${\displaystyle \epsilon_K=\frac{1}{\sqrt{2}}\sum_{|k|\geq
K+1}\mu_n(c)b_k(f)
\psi_{n,c}\left(\frac{k\pi}{c}\right).}$ Moreover,  from the well
known asymptotic behavior of the $\psi_{n,c}(x),$ for large values
of $x,$ see for example \cite{Karoui1}, one can easily check  that
${\displaystyle \epsilon_K =
o\left(\frac{1}{((K+1)\pi)^{1+s}}\right).}$ This computational
method of the $a_n(f)$  has the advantage to work for small as
well as large values of the smoothness coefficient $s>0.$\\

Also, note that if $f\in H^s([-1,1]),$ where $s > 1/2+ 2m, m\geq 1,$ is an
integer,  then $f\in C^{2m}([-1,1]).$ Moreover since
$\psi_{n,c}\in C^{\infty}(\mathbb R),$ then the classical Gaussian
quadrature method, see for example \cite{Andrews} gives us the
following approximate value $\widetilde a_n(f)$ of the $(n+1)-$th
expansion coefficient $a_n(f)=<f,\psi_{n,c}>,$
\begin{equation}\label{eq3.8}
\widetilde a_n(f)=\sum_{l=1}^m \omega_l f(x_l)\psi_{n,c}(x_l)=
a_n(f)+\epsilon_n,
\end{equation}
with  ${\displaystyle |\epsilon_n| \leq \sup_{\eta\in
[-1,1]}\frac{1}{b_m^2}\frac{(f \cdot
\psi_{n,c})^{(2m)}(\eta)}{(2m)!}.}$ Here, $b_m$ is the highest
coefficient of $\overline{P_m},$ and the different weights
$\omega_l$ and nodes $x_l,$ are easily computed by the special
method given in \cite{Andrews}.

The following examples illustrate the quality of approximation in $H^s(I)$ by the PSWFs.\\

\noindent {\bf Example 3:} In this example, we consider   the
Weierstrass function
\begin{equation}\label{eq5.1}
W_s(x)= \sum_{k\geq 0} \frac{\cos(2^{k}x)}{2^{ks}},\quad -1\leq
x\leq 1.
\end{equation}
Note that $W_s \in H^{s-\epsilon}([-1,1]),\,\forall \epsilon <
s,\, s
>0.$ We have considered the value of $c=100,$ and computed
$W_{s,N},$ the $N-$th terms truncated PSWFs series expansion of
$W_{s}$ with different values of ${\displaystyle \frac{3}{4}\leq s
\leq 2}$ and different values of $20\leq N \leq 100.$  Also, for
each pair $(s,N),$  we have computed the corresponding approximate
$L^2-$ error bound ${\displaystyle
E_N(s)=\left[\frac{1}{50}\sum_{k=-50}^{50}
(W_{s,N}(k/50)-W_{s}(k/50))^2\right]^{1/2}}.$  Table 2 lists the
obtained values of $E_N(s).$ Note that the numerical results given
by Table 2, follow what has been predicted by the theoretical
results of the previous section. In fact, the $L^2-$errors
$\|W_s-\Pi_N W_s\|_2$ is of order $O(N^{-s}),$ whenever
${\displaystyle N \geq N_c \sim \left[\frac{ 2c}{\pi}\right]+4.}$ In the case,
where $c=100,$ $N_c=67.$ The graphs of $W_{3/4}(x)$ and $W_{3/4,
N}(x),\, N=90$ are given by Figure 3.
\begin{center}
\begin{tiny}
{\tiny
\begin{table}[]
\caption{Values of $E_N(s)$ for various values of $N$ and $s.$}
\vskip 0.2cm
\begin{tabular}{ccccccc} \hline
    & $s=0.75$ &              $s= 1$    &     $s= 1.25$ &   $ s=1.5$ &   $ s=1.75$ & $s= 2.0$ \\ \hline
$N$ &  $E_N(s)$ & $E_N(s)$ &$E_N(s)$ &$E_N(s)$ &$E_N(s)$ &$E_N(s)$
\\ \hline
20& 4.57329E-01& 4.66173E-01 & 4.85990E-01  &    5.05973E-01  &    5.23232E-01  &  5.37227E-01 \\
30& 3.15869E-01& 3.11677E-01 & 3.28241E-01  &    3.48562E-01  &    3.67260E-01  &  3.82963E-01 \\
40& 1.06843E-01& 1.52009E-01 & 1.91237E-01  &    2.20969E-01  &    2.43432E-01  &  2.60523E-01 \\
50& 4.09844E-02& 6.88472E-02 & 1.01827E-01  &    1.26518E-01  &    1.44809E-01  &  1.58520E-01 \\
60& 3.30178E-02& 2.09084E-02 & 3.25551E-02  &    4.28999E-02  &    5.06959E-02  &  5.65531E-02  \\
70& 3.15097E-02& 8.82446E-03 & 2.51157E-03  &    7.35725E-04  &    2.33066E-04  &  1.04137E-04 \\
80& 3.01566E-02& 8.55598E-03 & 2.40312E-03  &    6.87458E-04  &    1.98993E-04  &  5.80481E-05  \\
90& 2.67972E-02& 7.64167E-03 & 2.14661E-03  &    6.15062E-04  &    1.78461E-04  &  5.22848E-05 \\
100&2.39141E-02& 6.72825E-03 & 1.82818E-03  &    5.10057E-04  &    1.45036E-04  &  4.19238E-05  \\

 \hline
\end{tabular}\end{table}
}
\end{tiny}
\end{center}

\begin{figure}[h]\hspace*{0.5cm}
{\includegraphics[width=14.5cm,height=5.5cm]{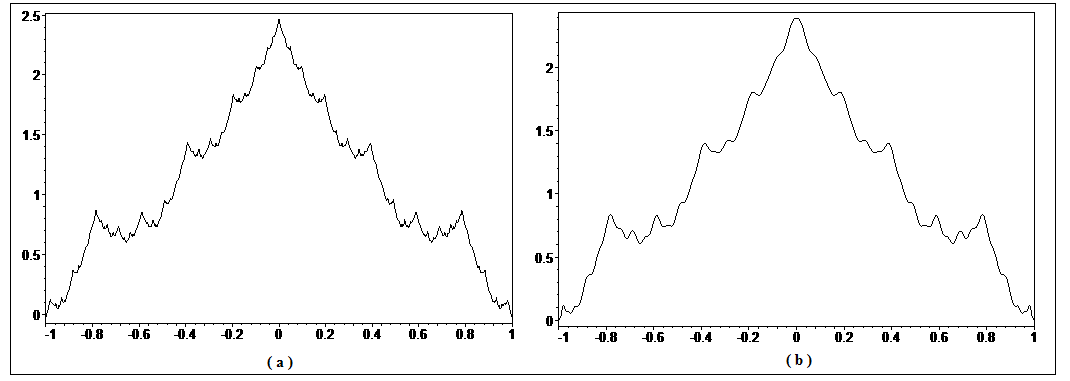}} \vskip
-0.5cm\hspace*{2cm}\caption{(a) graph of $W_{3/4}(x),\quad$ (b)
graph of $W_{3/4,N}(x), N=90.$ }
\end{figure}

\noindent {\bf Example 4:}  In this example, we let $s>0$ be any
positive real number and we consider the Brownian motion  function
$B_s(x)$ given by as follows.
\begin{equation}\label{eq5.1}
B_s(x)= \sum_{k\geq 1} \frac{X_k}{k^s} \cos(k\pi x),\quad -1\leq
x\leq 1.
\end{equation}
Here, $X_k$ is a Gaussian random variable. It is well known that
$B_s \in H^s([-1,1]).$ For the special case $s=1,$  we consider
the band-width $c=100,$ a truncation order $N=80$  and compute
$B_{1,N}$ the  approximation of $B_1$ by its $N-$th terms
truncated PSWFs series expansion. The graphs of $B_1$ and
$B_{1,N}$ are given by Figure 4.

\begin{figure}[h]\hspace*{0.5cm}
{\includegraphics[width=14.5cm,height=6.2cm]{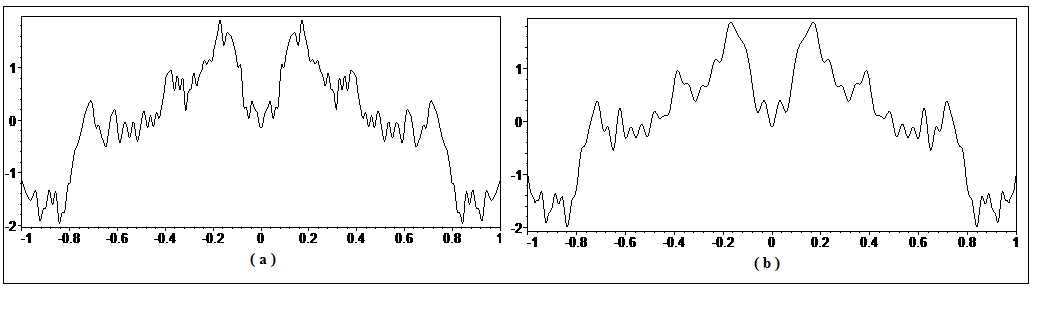}} \vskip
-0.5cm\hspace*{2cm}\caption{(a) graph of $B_1(x),\quad$ (b) graph
of $B_{1,N}(x), N=80.$ }
\end{figure}


\begin{remark}
From the quality of approximation in the Sobolev spaces
$H^s([-1,1])$ given in this paper and in \cite{Boyd1, Chen, Wang},
one concludes that for any value of the bandwidth $c\geq 0,$  the
approximation error $\|f-S_N f\|_2$ has the asymptotic order
$O(N^{-s}).$ Nonetheless, for a given $f\in H^s([-1,1]), s>0$
which we may assume to have a unit $L^2-$norm and for a given
error tolerance $\epsilon,$ the appropriate value of the bandwidth
$c\geq 0,$ corresponding to the minimum truncation order $N,$
ensuring that $\|f-S_N f\|_2\leq \epsilon,$ depends on whether or
not, $f$ has some significant Fourier expansion coefficients,
corresponding to large frequency components. In other words, the
faster decay to zero  of the Fourier coefficients of $f,$ the
smaller the value of the bandwidth should be and vice versa.
\end{remark}


\begin{thebibliography}{9999}

\bibitem{Abramowitz} M. Abramowitz and I. A. Stegun, Handbook of mathematical
 functions, Dover Publication, INC, New York.1972. pp.773-792.

\bibitem{Andrews} G. E. Andrews, R. Askey and  R. Roy, Special Functions,
Cambridge University Press, 2000.

\bibitem{Beylkin} G. Beylkin and L. Monzon, On Generalized Gaussian Quadrature for Exponentiels and their
Applications, {\it Appl. Comput. Harmon. Anal.} {\bf 12}, (2002),
332--373.

\bibitem{Bonami-Karoui1} A. Bonami and A. Karoui,  Useful bounds and eigenvalues decay  of the prolate spheroidal wave functions,
{\it  C. R. Math. Acad. Sci. Paris. Ser. I,} {\bf 352}  (2014), 229--234.

\bibitem{Bonami-Karoui2} A. Bonami and A. Karoui, Uniform Estimates of the Prolate
Spheroidal Wave Functions, submitted for publication (2014), available at http://arxiv.org/abs/1405.3676

\bibitem{Boyd1} J. P. Boyd, Prolate spheroidal wave functions as an alternative to Chebyshev and Legendre
 polynomials for spectral element and pseudo-spectral algorithms, {\it J. Comput. Phys.} {\bf 199}, (2004), 688--716.

\bibitem{Boyd2} J. P. Boyd, Approximation of an analytic function on a finite real interval by a bandlimited function
and conjectures on properties of prolate spheroidal functions,
{\it Appl. Comput. Harmon. Anal.} {\bf 25}, No.2,  (2003),
168--176.

\bibitem{Chen} Q. Chen, D. Gottlieb and J. S. Hesthaven, Spectral methods based on prolate spheroidal wave functions
for hyperbolic PDEs, {\it SIAM J. Numer. Anal.,} {\bf 43}, No. 5,
(2005), pp. 1912--1933.



\bibitem{Flammer} C. Flammer, {\it Spheroidal Wave Functions,} Stanford Univ. Press, CA, 1957.

\bibitem{Goss}  L. Gosse,  Compressed sensing with preconditioning for sparse recovery with subsampled matrices of Slepian prolate functions. {\sl Ann. Univ. Ferrara Sez. VII Sci. Mat.}, {\bf 59} (2013), 81--116.

\bibitem{Jamesson} G. J. O. Jameson, Elliptic integrals, the arithmetic-geometric mean and the
Brent-Salamin algorithm for $\pi,$ Notes, Dept. of Mathematics and Statistics, Lancaster University, Lancaster, U.K.

\bibitem{Karoui1} A. Karoui and T. Moumni, New efficient
methods of computing the prolate spheroidal wave functions and
their corresponding eigenvalues, {\it Appl. Comput. Harmon. Anal.}
{\bf 24}, No.3,  (2008), 269--289.

\bibitem{Landau1} H. J. Landau, The eigenvalue behavior of certain convolution equations,
{\it Trans. Amer. Math. Soc.,} {\bf 115,} (1965), 242--256.

\bibitem{Landau} H. J. Landau and H. O. Pollak, Prolate spheroidal  wave functions, Fourier analysis and
uncertainty-III. The dimension of space of essentially time-and
band-limited signals, {\it Bell System Tech. J.} {\bf 41}, (1962),
1295--1336.

\bibitem{LW} H. J. Landau and H. Widom, Eigenvalue
distribution of time and frequency limiting, {\it J. Math. Anal.
Appl.,} {\bf 77,} (1980), 469--481.

\bibitem{Li} L. W. Li, X. K. Kang, M. S. Leong, Spheroidal wave functions in electromagnetic
theory, Wiley-Interscience publication, 2001.


\bibitem{Lin} W. Lin, N. Kovvali and L. Carin, Pseudospectral method based on prolate spheroidal wave functions for
semiconductor nanodevice simulation, {\it Computer Physics
Communications,} {\bf 175} (2006), pp. 78--85.

\bibitem{Logan} J. A. Logan and J. D. Lakey, {\it Duration and Bandwidth Limiting: Prolate Functions, Sampling, and Applications,}
Applied and Numerical Harmonic Analysis Series, Birkh\"aser, Springer, New York, London, 2013.

\bibitem{Moore} I. C. Moore and M. Cada, Prolate spheroidal wave functions, an introduction to the Slepian
series and its properties, {\it Appl. Comput. Harmon. Anal.} {\bf
16}, No.3,  (2004), 208--230.

\bibitem{Nikoforov} A. N. Nikoforov and V. B. Uvarov, Special
functions of mathematical physics, translated from the Russian
edition, Birkh\"aser Verlag Basel, (1988).

\bibitem{Niven} C. Niven, On the Conduction of Heat in Ellipsoids of
Revolution, {\it Phil. Trans. R. Soc. Lond.,}  {\bf 171,} (1880),
117-151.

\bibitem{Osipov} A. Osipov, Certain inequalities involving prolate spheroidal wave functions and associated quantities,
{\it Appl. Comput. Harmon. Anal.}, {\bf 35},  (2013), 359--393.

\bibitem{Osipov2} A. Osipov, Certain upper bounds on the eigenvalues associated with
prolate spheroidal wave functions,
{\it Appl. Comput. Harmon. Anal.}, {\bf 35},  (2013), 309--340.


\bibitem{Osipov3} A. Osipov, V. Rokhlin and H. Xiao, {\it Prolate spheroidal wave functions of order zero. Mathematical tools for bandlimited approximation,}  Applied Mathematical Sciences, {\bf 187,} Springer, New York, 2013. 


\bibitem{Rokhlin} V. Rokhlin and H. Xiao, Approximate formulae
for certain prolate spheroidal wave functions valid for large
values of both order and band-limit, {\it Appl. Comput. Harmon.
Anal.} {\bf 22},  (2007), 105--123.

\bibitem{Rokhlin2} Y. Shkolnisky, M. Tygert and V. Rokhlin,
Approximation of bandlimited functions, {\it Appl. Comput. Harm.
Anal.,} {\bf 21}, (3), (2006), 413--420.


\bibitem{Slepian1} D. Slepian and H. O. Pollak, Prolate spheroidal wave functions, Fourier analysis and
uncertainty I, {\it Bell System Tech. J.} {\bf 40} (1961), 43--64.

\bibitem{Slepian3} D. Slepian,  Prolate spheroidal wave functions, Fourier analysis and
uncertainty--IV: Extensions to many dimensions; generalized
prolate spheroidal functions, {\it Bell System Tech. J.} {\bf 43}
(1964), 3009--3057.

\bibitem{Slepian2} D. Slepian, Some Asymptotic Expansions for Prolate Spheroidal Wave Functions,
{\it J. Math. Phys.,} {\bf 44}, No. 2, (1965), 99--140.


\bibitem{Wang} L. L. Wang,  Analysis of spectral approximations using prolate spheroidal wave functions.
 Math. Comp. 79 (2010), no. 270, 807--827.

 \bibitem{Wang2} L. L. Wang and J. Zhang, A new generalization of the PSWFs with applications to spectral
 approximations on quasi-uniform grids, {\it Appl. Comput. Harmon. Anal.}
 {\bf 29},  (2010), 303--329.

 \bibitem{Widom} H. Widom, Asymptotic behavior of the eigenvalues of certain integral equations. II. {\it Arc. Rational Mech. Anal.,}
  {\bf 17}  (1964), 215--229.\\

\bibitem{Xiao} H. Xiao, V. Rokhlin and N. Yarvin, Prolate spheroidal wave functions, quadrature and
interpolation, {\it Inverse Problems,} {\bf 17}, (2001), 805--838.






\end{thebibliography}
\end{document}